\documentclass[letterpaper, 9 pt, conference]{ieeeconf}      

\IEEEoverridecommandlockouts                              

\usepackage{graphics} 
\usepackage{epsfig} 
\usepackage{times} 
\usepackage{amsmath} 
\usepackage{amssymb}  
\usepackage{subfigure}
\usepackage{url}
\usepackage{algorithm}
\usepackage{dsfont}
\newtheorem{theorem}{Theorem}

\usepackage{cite}

\usepackage{color}
\usepackage{algorithm}
\newtheorem{properties}[theorem]{Properties}

\newtheorem{definition}[theorem]{Definition}

{ 
\newtheorem{remark}[theorem]{Remark}
\newtheorem{example}[theorem]{Example}
}

\newcommand{\bi}{\begin{itemize}}
\newcommand{\ei}{\end{itemize}}
\newcommand{\bd}{\begin{displaymath}}
\newcommand{\ed}{\end{displaymath}}
\newcommand{\be}{\begin{eqnarray*}}
\newcommand{\ee}{\end{eqnarray*}}

\title{\LARGE \bf
On Data-Driven Computation of Information Transfer for Causal Inference in Dynamical Systems
}
\author{Subhrajit Sinha and  Umesh Vaidya\\
\thanks{Financial support from the National Science Foundation grant  ECCS-1150405 and CNS-1329915 is gratefully acknowledged. S. Sinha and U. Vaidya is with the Department of Electrical \& Computer Engineering,
Iowa State University, Ames, IA 50011
        {\tt\small ugvaidya@iastate.edu}}%
}

\begin{document}
\maketitle
\begin{abstract}
In this paper, we provide a novel approach to capture causal interaction in a dynamical system from time-series data. In \cite{sinha_IT_CDC2016}, we have shown that the existing measures of information transfer, namely directed information, granger causality and transfer entropy fail to capture true causal interaction in dynamical system and proposed a new definition of information transfer that captures true causal interaction. The main contribution of this paper is to show that the proposed definition of information transfer in \cite{sinha_IT_CDC2016}\cite{sinha_IT_ICC} can be computed from time-series data. We use transfer operator theoretic framework involving Perron-Frobenius and Koopman operators for the data-driven approximation of the system dynamics and for the computation of information transfer. Several examples involving linear and nonlinear system dynamics are presented to verify the efficiency of the developed algorithm.

\end{abstract}

\section{Introduction}
Causality and influence characterization is an important problem in many different disciplines like economics, biological networks, social media, finance etc. Studying the cause and effects in these networks allows one to identify the influential and redundant nodes and thus aid immensely in the analysis of these systems. However, in many of these applications, the underlying mathematical model of the system is not available and hence one has to resort to causal inference from time series data. This is challenging even in bivariate case \cite{mooij2016distinguishing} and identification of causal structure from time series data is an active area of research. Two of the common methods for studying the cause-effect relationship in networks are graphical model approach \cite{pearl_book,sprites_book_causation} and system theoretic approach \cite{salapaka_reconstruction,adebayo,chow_sparse}. In another approach, concepts of information theory are used in such applications and a study of the information flow between the components of the network throws light on causality and the influential nodes of the network. In \cite{IT_socialmedia} the authors use information based metric to characterize the most influential nodes in social networks. In neuroscience, concepts of information theory are used to understand how information flows in different parts of the brain \cite{IT_brain} and identifying influence in gene regulatory networks \cite{IT_bionetwork1,IT_bionetwork2} . In economic and financial networks, information transfer can be used to infer causal interactions from the time series data \cite{granger_economics, IT_economics,granger_causality}.
Causality characterization was initially geared towards time series data and Granger causality \cite{granger_economics},\cite{granger_causality}, directed information \cite{IT_massey_directed},\cite{IT_kramer_directedit} and Schreiber's transfer entropy \cite{IT_schreiber} have been the most popular tools used for inferring the causality structure and influence characterization \cite{tatonetti2012data,kuhnert2014data,shu2013data,zou2009granger}. For details see \cite{mooij2016distinguishing} and references therein. However, in \cite{sinha_IT_CDC2016}, it was pointed out that all the above mentioned measures suffers from serious drawbacks and they fail to capture the correct causal structure even in very simple linear systems. Moreover, the authors also provided a new definition of information transfer in dynamical systems \cite{sinha_IT_CDC2016},\cite{sinha_IT_ICC} and had shown that this measure does capture the correct causal structure. 

The main idea behind the definition of information transfer proposed in \cite{liang_kleeman_prl} is the concept of freezing part of the system dynamics and has similarity with the definition of information transfer proposed in \cite{sinha_IT_ICC,sinha_IT_CDC2016}. The concept of freezing elevates the problem associated with other information-based causality measure and  captures the true causal structure \cite{sinha_IT_CDC2016}. However,  the freezing concept poses a challenge to compute our proposed information transfer from time series data. This is because the act of freezing the dynamics is analogous to intervening or modifying the system and it seems difficult to do once we have time-series data from the original un-freezed system. The second main challenge is inferring causal interaction from time-series data is the presence of noise in the data.

The main contribution of this paper is to provide data-driven approach for computing information transfer for inferring  causal interaction in dynamical system in noisy environment.  We use transfer operator theoretic framework involving transfer Perron-Frobenius ((P-F) and Koopman operators for the data-driven computation of information transfer. The transfer operators provides linear representation of the nonlinear systems by shifting focus from state space to the space of functions and measures. This linear representation is advantageous for the approximation of system dynamics from the time-series data. More recently there has been spur of research activities in the data-driven approximation of transfer operator. In particular, Dynamic Mode Decomposition (DMD), Extended Dynamic Mode Decomposition (EDMD), and Naturally Structured DMD (NSDDM) are some of the algorithms that are proposed for the finite dimensional approximation of these operators. In this paper we propose a robust implementation of NSDMD algorithm for the finite dimensional approximation of P-F and Koopman operators. The robust implementation is specially motivated to address the problem associated with the presence of noise in the data set.

  

The paper is organised as follows. In Section \ref{section_IT}, we provide the main definition of information transfer in dynamical system as developed in \cite{sinha_IT_ICC,sinha_IT_CDC2016}. In Section \ref{section_opt}, we provide a overview of main results from \cite{robust_dmd_acc,umesh_website} on robust approximation of transfer Koopman and P-F operator from time-series data. The robust approximation of transfer operators are used for the computation of information transfer and causality inference for linear system in Section \ref{section_linear}. Section \ref{section_general} discuss the causality inference results for nonlinear system. Finally we conclude the paper in Section \ref{section_conclusion}. 


\section{Information Transfer in Dynamical Systems}\label{section_IT}

In this section, we review the basics of information transfer in a dynamical system. For details, we refer the reader to \cite{sinha_IT_CDC2016, sinha_IT_ICC}. 
Consider the dynamical system $z(t+1)=F(z(t))+\xi(t)$, where $F = [F_x^\top \quad F_y^\top]^\top$, such that
{\small
\begin{eqnarray}\label{system2d}\left.
\begin{array}{ccl}
x(t+1) &=& F_x(x(t),y(t))+\xi_x(t)\\
y(t+1) &=&F_y(x(t),y(t))+\xi_y(t)
\end{array}\right\}
\end{eqnarray}
}
where $x\in\mathbb{R}^{|x|}$, $y\in\mathbb{R}^{|y|}$ (here $|\cdot|$ denotes the dimension of $\{\cdot\}$), $z=(x^\top,y^\top)^\top$, and  $F_x : \mathbb{R}^{|x|+|y|}\to\mathbb{R}^{|x|}$, $F_y : \mathbb{R}^{|x|+|y|}\to\mathbb{R}^{|y|}$ are assumed to be continuously differentiable and $\xi(t) = (\xi_x(t)^\top,\xi_y(t)^\top)^\top$ is additive independent and identically distributed noise. With some abuse of notation, we denote by $\rho(z(0))$ probability density function of $z$ at initial time $0$, which represents the uncertainty associated with initial condition $z(0)$. Similarly, we denote by $\rho(y(t+1)|y(t))$ the conditional probability density function of $y$ time $t+1$ conditioned on distribution of $y$ at time $t$. Next consider the following dynamical system where the dynamics in $x$ coordinate is frozen, going from time step $t$ to $t+1$ and defined as follows:
\begin{eqnarray}\label{system_xfreeze}\left.
\begin{array}{ccl}
x(t+1) &=& x(t)\\ 
y(t+1) &=& F_y(x(t),y(t))  + \xi_y(t)
\end{array}\right\}=F_{\not x}(z(t)) +\xi_{\not{x}}(t)
\end{eqnarray}
We denote by $\rho_{\not x}(y(t+1)|y(t))$ the probability density function of $y(t+1)$ conditional on $y(t)$ with the dynamics in $x$ coordinate freezed in time going from time step $t$ to $t+1$ as in Eq. (\ref{system_xfreeze}). We have following definition of information transfer from $x\to y$ going from time step $t$ to $t+1$. 
\begin{definition}\label{IT_def}[Information transfer] \cite{sinha_IT_CDC2016,sinha_IT_ICC} The information transfer from $x$ to $y$ for the dynamical system (\ref{system2d}), as the system evolves from time $t$ to time $t+1$ (denoted by $[T_{x\to y}]$), is given by following formula
\begin{small}
\begin{eqnarray}
[T_{x\to y}]=H(\rho(y(t+1)|y(t)))-H(\rho_{\not{x}}(y(t+1)|y(t))\label{IT}
\end{eqnarray}
\end{small}
where $H(\rho(y))=- \int_{\mathbb{R}^{|y|}} \rho(y)\log \rho(y)dy$ is the entropy of probability density function $\rho(y)$ and $H(\rho_{\not{x}}(y(t+1)|y(t))$ is the entropy of $y(t+1)$, conditioned on $y(t)$, where $x$ has been frozen as in Eq. (\ref{system_xfreeze}). 
\end{definition}
Using the relation between joint entropy and conditional entropy and using the fact that $H(y(t))=H_{\not{x}}(y(t))$, the information transfer from $x$ to $y$ can be written as
\begin{eqnarray}\label{info_trnasfer_eqv}
T_{x\to y} = H(y(t+1),y(t)) - H_{\not{x}}(y(t+1),y(t))
\end{eqnarray}
The above definition of information transfer can be extended to system with three subspace. In particular, let $x=(x_1^\top,x_2^\top)^\top$. We have following definition of information transfer in three subspace case. 
\begin{definition}
The information transfer from $x_1$ subspace to $y$, as the system evolves from time step $t$ to time step $t+1$ is 
\begin{small}
\begin{eqnarray*}
T_{x_1\to y} = H(y(t+1)|y(t)) - H_{\not{x}_1}(y(t+1)|y(t))
\end{eqnarray*}
\end{small}
\end{definition}

The information transfer from $x$ to $y$ depicts how evolution of $x$ affects the evolution of $y$, that is, it gives a quantitative measurement of the influence of $x$ on $y$. In lieu with this, we say that $x$ causes $y$ or $x$ influence $y$ if and only if the information transfer from $x$ to $y$ is non-zero and thus we have the following definition of influence in a dynamical system.

\begin{definition}\label{def_influence}
A state (or subspace) $x$ influences a state (or subspace) $y$ if and only if the information transfer from $x$ to $y$ is non-zero. 
\end{definition}
The objective of this paper is to discover the data-driven approach for inferencing causal interaction in dynamical system. Using the definition of influence, we will make use of information transfer measure to infer causality and influence. From the information transfer formula, $T_{x\to y}$, in Eq. (\ref{IT}), we notice that to compute the information transfer we are required to know the evolution of conditional probability density function, $\rho(z(t+1)|z(t))$, under two different dynamical systems i.e., the original dynamical system in Eq. (\ref{system2d}) and dynamical system with the dynamics in $x$ coordinates is freezed in Eq. (\ref{system_xfreeze}). The propagation of probability density function under system dynamics is governed by the linear transfer Perron-Frobenius operator \cite{Lasota}. However, for linear systems, with the assumption that the additive noise is i.i.d. Gaussian, one can derive an analytic expression for the information transfer between the states of the system \cite{sinha_IT_CDC2016}\cite{sinha_IT_ICC}.  

Consider the following linear time invariant dynamical system 
\begin{eqnarray}
z(t+1)=Az(t)+\sigma \xi(t)\label{lti}
\end{eqnarray}
where $z(t)\in \mathbb{R}^N$ and $\xi(t)$ is vector valued Gaussian random variable with zero mean and unit variance. 
We assume that the initial conditions have Gaussian distribution with covariance $\Sigma(0)$. Since the system is linear, the distribution of the system state for all future time will remain Gaussian with covariance $\Sigma(t)$ satisfying 
\begin{eqnarray}\label{evol_covar}
A \Sigma(t-1)A^\top+\sigma^2 I=\Sigma(t)
\end{eqnarray}
In particular, the distribution at any time $t$ is given by 
\begin{eqnarray*}
&&\rho_t(z) = \frac{1}{\sqrt{(2\pi)^N|\Sigma(t)|}}\times\\
&&\quad \exp [-\frac{1}{2}(z(t)-\mu(t))^\top \Sigma(t)^{-1}(z(t)-\mu(t))]
\end{eqnarray*}
where $\mu(t)$ is the mean at time $t$ and $\mu(t)=A\mu(t-1)$. This is because we have assumed that the additive Gaussian noise is zero mean i.i.d. noise.

To define the information transfer between various subspace we introduce following notation to split the $A$ matrix :
\begin{eqnarray}
z(t+1)=\begin{pmatrix}x^{'}\\y^{'}\end{pmatrix}=\begin{pmatrix}A_x&A_{xy}\\ A_{yx}&A_{y}\end{pmatrix}\begin{pmatrix}x\\y\end{pmatrix}+\sigma \xi\label{splittingxy}
\end{eqnarray}

The $A$ matrix can be further split using the subspace decomposition $x=(x_1^\top,x_2^\top)^\top$ as follows:
\begin{eqnarray}
\begin{pmatrix}A_x&A_{xy}\\ A_{yx}&A_{y}\end{pmatrix}=\begin{pmatrix}A_{x_1}&A_{x_1x_2}& A_{x_1 y}\\A_{x_2x_1}&A_{x_2}& A_{x_2 y}\\ A_{y x_1}&A_{y x_2}& A_{y}\end{pmatrix}\label{splittingA}
\end{eqnarray}
Based on the decomposition of the system $A$ matrix we can also decompose the covariance matrix $\Sigma$ at time instant $t$ as follows. 
\begin{eqnarray}
\Sigma=\begin{pmatrix}\Sigma_x&\Sigma_{xy}\\\Sigma_{xy}^\top& \Sigma_y\end{pmatrix}=\begin{pmatrix} \Sigma_{x_1}&\Sigma_{x_1x_2}&\Sigma_{x_1 y}\\\Sigma_{x_1x_2}^\top&\Sigma_{x_2}&\Sigma_{x_2 y}\\\Sigma_{x_1y}^\top&\Sigma_{x_2y}^\top&\Sigma_{y}\end{pmatrix}\nonumber\\
\label{sigma_dec}
\end{eqnarray}
Using the above notation, we state following theorem providing explicit expression for information transfer in linear dynamical system during transient and steady state.

\begin{theorem}\cite{sinha_IT_CDC2016}\cite{sinha_IT_ICC}.
Consider the linear dynamical system (\ref{lti}) and associated splitting of state space in Eqs. (\ref{splittingxy}) and (\ref{splittingA}). We have following expression for information transfer between various subspace
\begin{eqnarray}
[T_{x\to y}]_t^{t+1}=\frac{1}{2}\log \frac{|A_{yx}\Sigma^s_{y}(t)A_{yx}^\top +\sigma^2 I_y|}{|\sigma^2 I_y|}\label{transferxy}
\end{eqnarray}
where \begin{eqnarray}\Sigma^s_y(t)=\Sigma_x(t)-\Sigma_{xy}(t)\Sigma_y(t)^{-1}\Sigma_{xy}(t)^\top\label{schur_complement}
\end{eqnarray} is the Schur complement of $\Sigma_{y}(t)$ in the matrix $\Sigma(t)$ and $I_y$ is the identity matrix of dimension equal to dimension of $y$. 

\begin{eqnarray}
[T_{x_1\to y}]_t^{t+1}=\frac{1}{2}\log \frac{|A_{yx}\Sigma^s_y(t)A_{yx}^\top +\sigma^2 I_y |}{|A_{yx_2}(\Sigma_y^{s})_{yx_2}(t)A_{yx_2}^\top+\sigma^2 I_y|}\label{transferx1y}
\end{eqnarray}

where $|\cdot|$ is the determinant and $ (\Sigma_y^s)_{yx_2}$ is the Schur complement of $\Sigma_{y}$ in the matrix 
\[\begin{pmatrix}\Sigma_{x_2}&\Sigma_{x_2y}\\\Sigma_{x_2 y}^\top&\Sigma_y\end{pmatrix}\]
\end{theorem}
\[\]
\begin{remark}\label{remark_topology}
From the analytical expression of information transfer Eq. (\ref{transferx1y}) it is easy to show that the information transfer from state $T_{x_1\to y}$ is zero if and only if $A_{yx_1}$ is zero. This connection between the information transfer and the structure of the $A$ matrix can be used for network topology identification. We will discuss this results in the simulation section.  
\end{remark}

\section{Robust approximation of P-F and Koopman operators}\label{section_opt}
For the computation of information transfer we are required to know the evolution of probability density function. The evolution of probability density function in forward time is governed by transfer Perron-Frobenius (P-F) operator.
In this section we describe in brief the main results from \cite{robust_dmd_acc,umesh_website} on robust approximation of transfer Koopman and P-F operators from the time series data. Transfer operator theoretic framework involving Koopman and Perron-Frobenius (P-F) operators are linear operators \cite{Lasota}. 
Various methods are proposed for finite dimensional approximation of these infinite dimensional operators from time-series data. The most popular among them are Dynamic Mode Decomposition (DMD) and Extended DMD (EDMD) \cite{williams2015data,schmid2010dynamic}. However, these methods and others assume that the time series data is noise free and hence does not explicitly account for uncertainty in data-set. Here we propose a novel approach based on the theory of robust optimization for robust approximation of the transfer Koopman operator from noisy time series data \cite{robust_dmd_acc,umesh_website}. Consider a discrete-time dynamical system forced with stochastic input. 
\begin{eqnarray}
z_{t+1}=T(z_t,\xi_t)\label{system}
\end{eqnarray}
where $T:Z\times W \to  Z$ with $X\subset \mathbb{R}^N$ is assumed to be invertible with respect to $z$ for each fixed value of $\xi$ and smooth diffeomorphism. $\xi_t\in W$ is assumed to be independent identically distributed (i.i.d) random variable drawn from probability distribution $\vartheta$ i.e., 
\[{\rm Prob}(\xi_t\in B)=\vartheta(B)\]
for every set $B\subset W$ and all $t$.
Furthermore, we denote by ${\cal B}(Z)$ the Borel-$\sigma$ algebra on $X$ and ${\cal M}(Z)$ the vector space of bounded complex-valued measure on $X$.  Associated with this discrete time dynamical system are two linear operators namely Koopman and Perron-Frobenius (P-F) operator. These two operators are defined as follows.
\begin{definition}[Perron-Frobenius Operator]  $\mathbb{P}:{\cal M}(Z)\to {\cal M}(Z)$ is given by
{\small
\begin{eqnarray}
[\mathbb{P}\mu](A)=\int_{{\cal X} }\int_W\chi_{A}(T(x,v))d\vartheta(v)d\mu(x)=\int_X p(x,A)d\mu(x)
\end{eqnarray}
}
where $\chi_A(x)$ is the indicator function for set $A$ and $p(x,A)$ is the transition probability function.  \end{definition}
For deterministic dynamical system $p(x,A)=\delta_{T(x)}(A)$. Under the assumption that $p(x,\cdot)$ is absolutely continuous with respect to Lebesgue measure, $m$, we can write 
\[p(x,A)=\int_A p(x,y)dm(y)\]
for all $A\subset X$. Under this absolutely continuous assumption, the P-F operator on the space of densities $L_1(X)$ can be written as  \footnote{with some abuse of notation we are using the same notation for the P-F operator defined on the space of measure and densities.}
\[[\mathbb{P}g](y)=\int_X p(x,y)g(x)dm(x)\]
\begin{definition}[Invariant measures] Invariant measures are the fixed points of
the P-F operator $\mathbb{P}$ that are additionally probability measures. Let $\bar \mu$ be the invariant measure then, $\bar \mu$ satisfies
\[\mathbb{P}\bar \mu=\bar \mu\]
\end{definition}
Under the assumption that the state space $X$ is compact, it is known that the P-F operator admits at least one invariant measure.
\begin{definition} [Koopman Operator] Given any $h\in\cal{F}$, $\mathbb{U}:{\cal F}\to {\cal F}$ is defined by
\[[\mathbb{U} h](x)={\bf E}_{\xi}[h(T(x,\xi))]=\int_W h(T(x,v))d\vartheta(v)\]
\end{definition}

\begin{properties}\label{property}
Following properties for the Koopman and Perron-Frobenius operators can be stated.

\begin{enumerate}
\item [a).] For any function $h\in{\cal F}$ such that $h\geq 0$, we have $[\mathbb{U}h](x)\geq 0$ and hence Koopman is a positive operator.

\item [b).] If we define P-F operator act on the space of densities i.e., $L_1(X)$ and Koopman operator on space of $L_\infty(X)$ functions, then it can be shown that the P-F and Koopman operators are dual to each others as follows 
\begin{eqnarray*}
&&\left<\mathbb{U} f,g\right>=\left<f,\mathbb{P} g\right>
\end{eqnarray*}
where $f\in L_{\infty}(X)$ and $g\in L_1(X)$.

\item [c).] For $g(x)\geq 0$, $[\mathbb{P}g](x)\geq 0$.

\item [d).] Let $(X,{\cal B},\mu)$ be the measure space where $\mu$ is a positive but not necessarily the invariant measure, then the P-F operator  satisfies  following property.
\[\int_X [\mathbb{P}g](x)d\mu(x)=\int_X g(x)d\mu(x)\]\label{Markov_property}
\end{enumerate}
\end{properties}

We next discuss the approximation of these two operators from time series data. 
Consider snapshots of data set obtained from simulating a discrete time random dynamical system $z\to T(z,\xi)$ or from an experiment
\begin{eqnarray}
{\cal Z}  = [z_0,z_1,\ldots,z_M]
 \label{data}
\end{eqnarray}
where $z_i\in Z\subset \mathbb{R}^N$. The data-set $\{z_k\}$ can be viewed as sample path trajectory generated by random dynamical system and could be corrupted by either process or measurement noise or both.

A large number of sample path trajectories need to be simulated to realize sufficient statistics of the random dynamical system. However, in practice, only few sample path trajectories over finite time horizon are available, and it is hard to approximate the statistics of RDS using the limited amount of data-set. Furthermore, rarely one knows the probability distribution of the underlying noise process, i.e., $\vartheta$. Estimating $\vartheta$ is in itself a challenging problem. In spite of these difficulties, it is essential to develop an algorithm for the approximation of transfer operators that explicitly account for the uncertainty in data-set. We propose a robust optimization-based approach to address this challenge. In particular, we consider deterministic, but norm bounded uncertainty in the data set. Since the trajectory $\{z_k\}$ is one particular realization of the RDS, the other random realization can be assumed to be obtained by perturbing $\{z_k\}$. We assume that the data points $z_k$ are perturbed by norm bounded deterministic perturbation of the form 
\[\delta z_k=z_k+\delta,\;\;\; \delta\in \Delta.\]
Several possible choices for the uncertainty set $\Delta$ can be considered. For example
\[\Delta:=\{\delta \in \mathbb{R}^n:\;\; \parallel \delta\parallel_2\leq \rho\}\]
restrict the $2$-norm of $\delta$ to $\rho$. Another possible choice could be 
\[\Delta:=\{\delta \in \mathbb{R}^n:\;\; \parallel \delta\parallel_{Q_i}\leq 1,\;\;i=1,\ldots, d\}\]
where $Q_i\geq 0$ and implies that uncertainty $\delta$ lies at the intersection of ellipsoids. More generally, one can also consider $\Delta$ set to be of the form
\[\Delta=\{\delta\in \mathbb{R}^n: h_i(\delta)\leq 0,\;\;i=1,\ldots,d\}\]
for some convex function $h_i(\delta)$.

Now let $\mathcal{D}=
\{\psi_1,\psi_2,\ldots,\psi_K\}$ be the set of dictionary functions or observables. The dictionary functions are assumed to belong to $\psi_i\in L_2(X,{\cal B},\mu)={\cal G}$, where $\mu$ is some positive measure, not necessarily the invariant measure of $T$. Let ${\cal G}_{\cal D}$ denote the span of ${\cal D}$ such that ${\cal G}_{\cal D}\subset {\cal G}$. The choice of dictionary functions are very crucial and it should be rich enough to approximate the leading eigenfunctions of Koopman operator. Define vector valued function $\mathbf{\Psi}:X\to \mathbb{C}^{K}$ as
\begin{equation}
\mathbf{\Psi}(z):=\begin{bmatrix}\psi_1(z) & \psi_2(z) & \cdots & \psi_K(z)\end{bmatrix}\label{dic_function}
\end{equation}
In this application, $\mathbf{\Psi}$ is the mapping from physical space to feature space. Any function $\phi,\hat{\phi}\in \mathcal{G}_{\cal D}$ can be written as
\begin{eqnarray}
\phi = \sum_{k=1}^K a_k\psi_k=\boldsymbol{\Psi a},\quad \hat{\phi} = \sum_{k=1}^K \hat{a}_k\psi_k=\boldsymbol{\Psi \hat{a}}\label{expand}
\end{eqnarray}
for some set of coefficients $\boldsymbol{a},\boldsymbol{\hat{a}}\in \mathbb{C}^K$. 
Let \begin{eqnarray}
 \hat{\phi}(z)=[\mathbb{U}\phi](z)+r=E_\xi[\phi(T(z,\xi))]+r. \label{residual}
\end{eqnarray}
Unlike deterministic case where we evaluate (\ref{residual}) at the data point $\{z_k\}$, for the uncertain case  we do not have sufficient data points to evaluate the expected value in the above expression. Instead we use the fact that different realizations of the RDS will consist of the form $\{z_k+\delta\}$ with $\delta\in \Delta$ to write (\ref{residual}) as follows:
\begin{eqnarray}
\hat{\phi}(z_m + \delta z_m)=\phi(z_{m+1})+r, \;\;\;k=1,\ldots,M-1.
\end{eqnarray}
The objective is to minimize the residual for not just one pair of data points $\{z_m,z_{m+1}\}$, but over all possible pairs of data points of the form $\{z_m+\delta,z_{m+1}\}$. 
Using (\ref{expand}) we write the above as follows:
\[
\boldsymbol{\Psi}(z_k + \delta z_k)\boldsymbol {\hat{a}}=\boldsymbol{\Psi}(z_{k+1})\boldsymbol {a}+r.
\]

We seek to find matrix $\bf K$, the finite dimensional approximation of Koopman operator that maps coefficient vector $\boldsymbol{a}$ to $\boldsymbol{\hat{a}}$, i.e., ${\bf K}\boldsymbol{a}=\boldsymbol{\hat a}$, while minimizing the residual term, $r$.  Premultiplying by $\boldsymbol{\Psi}^\top( z_m)$ on both the sides of above expression and summing over $m$ we obtain 
\[\left[\frac{1}{M}\sum_m \boldsymbol{\Psi}^\top( z_m) \boldsymbol{\Psi}(z_m + \delta z_m){\bf K}-\boldsymbol{\Psi}^\top(z_m)\boldsymbol{\Psi}( z_{m+1})\right]{\boldsymbol a}.\]
In the absence of the uncertainty the objective is to minimize the appropriate norm of the  quantity inside the bracket over all possible choices of matrix $\bf K$. However, for robust approximation, presence of uncertainty acts as an adversary whose goal is to maximize the residual term. Hence the robust optimization problem can be formulated as a $\min-\max$ optimization problem as follows. 

\begin{equation}\label{edmd_robust}
\min\limits_{\bf K}\max_{\delta\in \Delta}\parallel {\bf G}_\delta{\bf K}-{\bf A}\parallel_F=:\min\limits_{\bf K}\max_{\delta\in \Delta} {\mathfrak{F}}({\bf K}, {\bf G}_\delta,{\bf A})
\end{equation}
where
\begin{eqnarray}\label{edmd1}
&&{\bf G}_\delta=\frac{1}{M}\sum_{m=1}^M \boldsymbol{\Psi}({ z}_m)^\top \boldsymbol{\Psi}({z_m + \delta z}_m)\nonumber\\
&&{\bf A}=\frac{1}{M}\sum_{m=1}^M \boldsymbol{\Psi}({ z}_m)^\top \boldsymbol{\Psi}({ z}_{m+1}),
\end{eqnarray}
with ${\bf K},{\bf G}_\delta,{\bf A}\in\mathbb{C}^{K\times K}$. The $\min-\max$ optimization problem (\ref{edmd_robust}) is in general nonconvex and will depend on the choice of dictionary functions. This is true because ${\mathfrak{F}}$ in (\ref{edmd_robust}) is not in general concave function of $\delta$ for fixed $\bf K$. 
Hence, we convexify the problem as follows

\begin{equation}\label{edmd_robust_convex}
\min\limits_{\bf K}\max_{\delta{\bf G}\in \bar \Delta}\parallel ({\bf G}+\delta {\bf G}){\bf K}-{\bf A}\parallel_F
\end{equation}
where $\delta {\bf G}\in \mathbb{R}^{K\times K}$ is the new perturbation term characterized by uncertainty set  $\bar \Delta$ which lies in the feature space of dictionary function and the matrix ${\bf G}=\frac{1}{M}\sum_{m=1}^M \boldsymbol{\Psi}({ z}_m)^\top \boldsymbol{\Psi}({x}_m)
$. $\bar \Delta$ is the new uncertainty set defined in the feature space and will inherit the structure from set $\Delta $ in the data space. In particular, it is easy to show that

\begin{eqnarray}
\parallel \delta G\parallel_F \leq \lambda \Lambda \Gamma
\end{eqnarray}
where $\parallel\delta z_m\parallel_F\leq \lambda$, $\parallel {\bf \Psi}(z_m)\parallel_F\leq \Lambda$ and $\parallel {\bf \Psi}'(z_m)\parallel_F\leq \Gamma$ for all $m$. 


In \cite{Umesh_NSDMD}, we proposed Naturally Structured Dynamic Mode Decomposition (NSDMD) algorithm for finite dimensional approximation of the transfer Koopman and P-F operator. Apart from preserving positivity and Markov properties of the transfer operator, this algorithm exploits the duality between P-F and Koopman operator to provide the approximation of P-F operator. 
The algorithm presented for the robust approximation of Koopman operator can be combined with NSDMD for the robust approximation of P-F operator. In particular, under the assumption that all the dictionary functions are positive, following modification can be made to optimization formulation (\ref{edmd_robust_convex}) for the approximation of Koopman operator. 

\begin{eqnarray}\label{edmd_robust_convexPF} \nonumber
&&\min\limits_{\bf K}\max_{\delta{\bf G}\in \bar \Delta}\parallel ({\bf G}+\delta {\bf G}){\bf K}-{\bf A}\parallel_F\nonumber\\ \nonumber
&&{\rm s.t.}\quad{\bf K}_{ij}\geq 0\\ \nonumber
&&\qquad [\Lambda {\bf K}\Lambda^{-1}]_{ij}\geq 0\\
&&\qquad \Lambda{\bf K}\Lambda^{-1}\mathds{1}=\mathds{1}
\end{eqnarray}
where $\Lambda = \langle\boldsymbol{\Psi}(z),\boldsymbol{\Psi}(z)\rangle$ with $[\Lambda]_{ij}=\langle\psi_i,\psi_j\rangle$ is symmetric positive definite  matrix.
We refer the interested reader to \cite{Umesh_NSDMD} for details of NSDMD formulation.
Using duality the  robust approximation of the P-F operator, $\bf P$, can then be written as ${\bf P}=\Lambda^{-1}{\bf K}^\top {\Lambda}$. Most common approach for solving the robust optimization problem is by using a robust counterpart. In the following  section we show that the robust counterpart of the robust optimization problem can be constructed and is a convex optimization problem. 

The robust optimization problem (\ref{edmd_robust_convex}) has some interesting connection with optimization problems involving regularization term. In particular, we have following Theorem.

\begin{theorem}
Following two optimization problems 
\begin{eqnarray}
\min\limits_{\bf K}\max_{{\delta {\bf G}:}\parallel \delta{\bf G}\parallel_F\leq \lambda}\parallel ({\bf G}+\delta {\bf G}){\bf K}-{\bf A}\parallel_F \label{ocp}
\end{eqnarray}
\begin{eqnarray}\label{rob_eqv}
\min\limits_{\bf K}\parallel {\bf G}{\bf K}-{\bf A}\parallel_F+\lambda \parallel {\bf K}\parallel_F\label{regular}
\end{eqnarray}
are equivalent.
\end{theorem}
Refer to \cite{caramanis201214} for the proof.

Using the above equivalence the robust implementation of NSDMD can be written as following optimization problem

\begin{eqnarray}
\min\limits_{\bf K}\parallel {\bf G}{\bf K}-{\bf A}\parallel_F+\lambda \parallel {\bf K}\parallel_F\label{cost}\\
{\rm s.t.}\left\{\begin{array}{ccl}\qquad{\bf K}_{ij}&\geq &0\\ 
\qquad[\Lambda {\bf K}\Lambda^{-1}]_{ij}&\geq & 0\\
\qquad\Lambda{\bf K}\Lambda^{-1}\mathds{1}&=&\mathds{1}\end{array}\right. \label{positive3}
\end{eqnarray}

Let $w(t)\in \mathbb{R}^K$ and $v(t)\in \mathbb{R}^K$ be a row vector and column vector respectively. Furthermore, $w_i(t)\geq 0$ and $\sum_{i}w_i=1$. We have,
{\small
\begin{eqnarray}
&&\textnormal{P-F Operator : } w(t+1) = w(t){\bf P}^\top\label{PFevolution}\\
&&\textnormal{Koopman Operator : } v(t+1) = {\bf K} v(t) \label{Koopmanevolution}
\end{eqnarray}}
The robust implementation of EDMD will simply correspond to minimizing the unconstrained cost function (\ref{regular}) without the positivity and Markov constraints (\ref{positive3}). In particular, the  optimization problem for robust EDMD can then be written as 
\begin{eqnarray}
\min\limits_{\bf K}\parallel {\bf G}{\bf K}-{\bf A}\parallel_F+\lambda \parallel {\bf K}\parallel_F\label{regular2}.
\end{eqnarray}

\section{Causal Inference in Linear Dynamical System}\label{section_linear}
In this section we outline the procedure for computing the information-based causal inference from time series data for the case of linear dynamical system. We will employ the data-driven approximation of transfer operators discussed in the previous section and also exploit the fact that the analytical expression for information transfer in linear system are available (i,e., Eqs. (\ref{transferxy}-\ref{transferx1y})).
For ease of understanding, we discuss the procedure of information transfer computation for a two dimensional or two subspace case, as given by system Eq. (\ref{splittingxy}). The general case will follow from the two subspace case.

Note that in all the subsequent discussion we assume that we have access to all the states of the system. The problem of computing the information transfer based on output measurements could be more realistic problem and is a topic of our ongoing investigation. Let the time series data be given by

\begin{eqnarray}
\mathcal{D}=\bigg[\begin{pmatrix}
x_0\\
y_0
\end{pmatrix},  \begin{pmatrix}
x_1\\
y_1
\end{pmatrix}, 
\cdots, \begin{pmatrix}
x_{M-1}\\
y_{M-1}
\end{pmatrix}\bigg] \label{data_original}
\end{eqnarray}

Since the data is assumed to be generated from a linear dynamical system, we use linear dictionary functions i.e., $\psi_k(z)=z_k$. Furthermore, 
the number of dictionary functions are taken to be  equal to the size of the system i.e., $N$. For the linear system case we use optimization formulation (\ref{regular2}) i.e., EDMD with linear dictionary functions for the approximation of Koopman operator. With the linear choice of dictionary function and number of dictionary function equal to size of the system it is not difficult to show that the approximation of Koopman operator, $\bf K$, is the system $A$ matrix itself. Let $\bar A=\bf K\in \mathbb{R}^{N\times N}$ be the estimated system dynamics obtained using optimization formulation (\ref{regular2}). Under the assumption that the initial covariance matrix is $\bar \Sigma(0)$, the propagation of the covariance matrix under the estimated system dynamics $\bar A$ is given by 
\begin{eqnarray}
\bar \Sigma(t)=\bar A\bar \Sigma(t-1)\bar {A}^\top+\sigma^2 I\label{covariance_prop}
\end{eqnarray}
Both $\bar A$ and $\bar \Sigma$ can be decomposed according to Eqs. (\ref{splittingA}) and (\ref{sigma_dec}). The conditional entropy $H(y_{t+1}|y_t)$ for the non-freeze case is computed using the following formula \cite{sinha_IT_CDC2016,sinha_IT_ICC}.  
\begin{eqnarray}\label{cond_entr}
H(y_{t+1}|y_t) = \frac{1}{2}\log |\bar A_{yx}\bar \Sigma_y^S(t)\bar A_{yx}^\top + \left(\frac{\lambda}{3}\right)^2 I|.
\end{eqnarray}
where $|\cdot |$ is the determinant, $\lambda$ is the bound on the process noise,  $\bar \Sigma_y^S(t)$ is the Schur complement of $y$ in the covariance matrix $\bar \Sigma(t)$ (refer to Eq. (\ref{schur_complement}) for Schur complement). In computing the entropy, we assume that the noise is i.i.d. Gaussian with covariance $\Sigma = \textnormal{diag}(\sigma^2,\cdots ,\sigma^2)$ so that one can take the bound as $\lambda = 3\sigma$, to cover the essential support of the Gaussian distribution. 

Computing the conditional entropy of $y$ when $x$ is frozen from the time series data obtained from the non-freeze dynamics is a challenge. To replicate the effect of $x$ freeze dynamics we modify the original data set (\ref{data_original}) as follows.  
\begin{small}
\begin{eqnarray}
\mathcal{D}_{\not{x}}=\bigg[\left\{\begin{pmatrix}
x_0\\
y_0
\end{pmatrix}, \begin{pmatrix}
x_0\\
y_1
\end{pmatrix}\right\}, \left\{\begin{pmatrix}
x_1\\
y_1
\end{pmatrix}, \begin{pmatrix}
x_1\\
y_2
\end{pmatrix}\right\},
\cdots, \nonumber\\
\cdots ,\left\{ \begin{pmatrix}
x_{M-1}\\
y_{M-1}
\end{pmatrix}, \begin{pmatrix}
x_{M-1}\\
y_M
\end{pmatrix}\right\}\bigg]\label{mod_data}
\end{eqnarray}
\end{small}
If the original data set has $M$ data points, then the modified data set has $(2M-2)$ data points. The idea is to find the best mapping that propagate points of the form $[x_{t-1} \quad y_{t-1}]^\top$ to $[x_{t-1} \quad y_t]^\top$  (i.e., $x$ freeze) for $t = 1, 2, \hdots , M$. The estimated dynamics $\bar A_{\not x}$, when $x$ is frozen,  is calculated using the optimization formulation (\ref{regular}) but this time  applied to the data set (\ref{mod_data}). Once the frozen model is calculated, the entropy $H_{\not{x}}(y_{t+1}|y_t)$ is calculated using exactly the same procedure outline for $H(y_{t+1}|y_t)$ but this time applied to $\bar A_{\not x}$. Finally the information transfer from $x\to y$ is computed using the formula
\[T_{x\to y}=H(y_{t+1}|y_t)-H_{\not{x}}(y_{t+1}|y_t)\]

The algorithm for computing the information transfer for the linear system case can be  summarized as follows:
\begin{algorithm}[htp!]
\caption{Information Transfer: Linear System}
\begin{enumerate}
\item{From the original data set (\ref{data_original}), compute the estimate of the system matrix $\bar A$ using the optimization formulation (\ref{regular2})}

\item{Assume $\bar \Sigma(0)$ and compute $\bar \Sigma(t)$ using Eq. (\ref{covariance_prop}). Determine $\bar A_{yx}$ and $\bar \Sigma_y^S$ to calculate the conditional entropy $H(y_{t+1}|y_t)$ using (\ref{cond_entr}).}
\item{From the original data set (\ref{data_original}) form the modified data set for the $x$ freeze dynamics as given by Eq. (\ref{mod_data}).}
\item{Follow steps (1)-(2) to compute the conditional entropy $H_{\not{x}}(y_{t+1}|y_t)$.}
\item{Compute the transfer $T_{x\to y}$ as $T_{x\to y} = H(y_{t+1}|y_t) - H_{\not{x}}(y_{t+1}|y_t)$.}
\end{enumerate}\label{algo_IT}
\end{algorithm}
\subsection{Simulation Results}\label{section_simulation}

\begin{figure}[htp!]
\centering
\includegraphics[scale=.33]{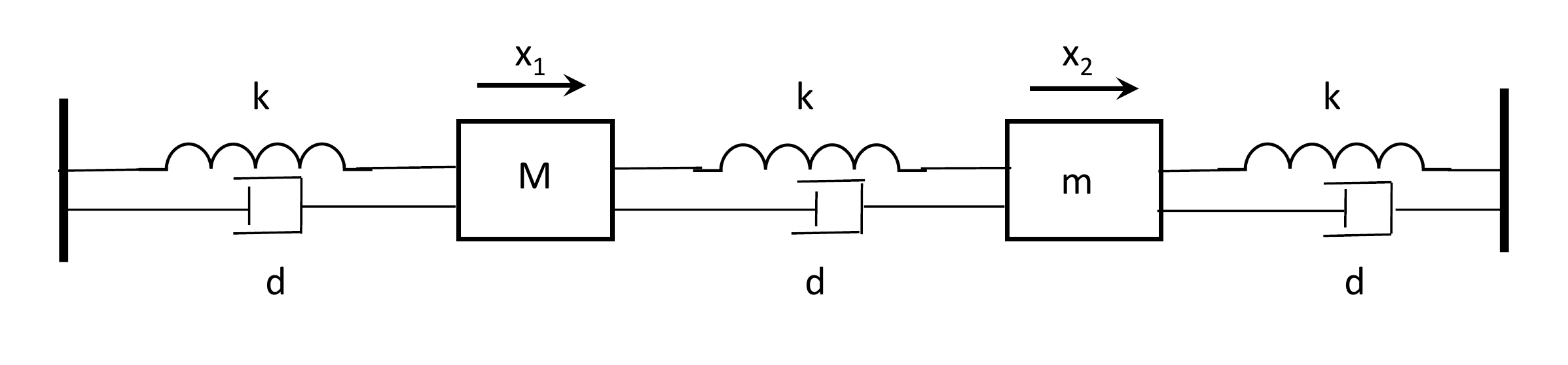}
\caption{Mass-spring-damper system}\label{mass_spring_fig}
\end{figure}

{\bf Example 1}: In the first example, we discuss the physical meaning of information transfer and demonstrate how it can be used to characterize influence in a dynamical system. Consider a mass-spring-damper system, as shown in Fig. \ref{mass_spring_fig}.

The equations of motion for the mass-spring system are 
\begin{eqnarray}
M\ddot{x}_1 + 2d\dot{x}_1 - d\dot{x}_2 + 2kx_1 - kx_2 = 0\\
m\ddot{x}_1 + 2d\dot{x}_2 - d\dot{x}_1 + 2kx_2 - kx_1 = 0
\end{eqnarray}
where $M,m$ are the masses, $d$ is the damping coefficient and $k$ is the spring constant. We assume that the damping coefficients of the dampers are equal and so are the spring constants of the springs.
In state space form, the system can be represented as
{\small
\begin{eqnarray*}
\begin{pmatrix}
\dot{z}_1\\
\dot{z}_2\\
\dot{z}_3\\
\dot{z}_4
\end{pmatrix} = \begin{pmatrix}
0 & 1 & 0 & 0\\
-2k/M & -2d/M & k/M & d/M\\
0 & 0 & 0 & 1\\
k/m & d/m & -2k/m & -2d/m
\end{pmatrix}\begin{pmatrix}
{z}_1\\
{z}_2\\
{z}_3\\
{z}_4
\end{pmatrix}
\end{eqnarray*}
}
where $z_1 = x_1$, $z_2 = \dot{x}_1$, $z_3 = x_2$ and $z_4=\dot{x}_2$.
For simulation purposes, we choose $M=10$, $m = 1$, $k=1$ and $d=5$ with appropriate units.  Since $M>m$, a perturbation (perturbed so that it has some non-zero initial velocity) in $M$ will result in larger oscillations in the masses, compared to the case when $m$ is perturbed by the same amount. Hence, we can conclude that $M$ has a large influence on $m$, whereas, $m$ has much smaller influence on $M$. In the language of information transfer between the states, this can be characterized by the information transfer from the position variable of one mass to the velocity variable of the other mass. From the analytical expression of information transfer, that is, equation (\ref{transferx1y}), the information transfer from $z_1\to z_4$ and $z_3\to z_2$ are $T_{z_1\to z_4} = 0.174$ and $T_{z_3\to z_2}=0.053$. This confirms the fact that $M$ has a much larger influence on $m$, whereas $m$ has negligible effect on $M$. 

For calculating these transfers from data, the system was initialized at $[1\quad .2\quad 1\quad .3]^\top$ and data was collected for 10 time steps with sampling time $\delta t = 0.1$ seconds. An additive Gaussian noise of variance 0.1 was added to the system. With this, the information transfer values were calculated as $T_{z_1\to z_4} = 0.1502$ and $T_{z_3\to z_2}=0.03$. Though the exact values do not match, but they are close and more importantly, they do convey the fact that $M$ has a much larger influence on $m$ and $m$ has negligible influence on $M$.

\begin{figure}[htp!]
\centering
\includegraphics[scale=.3]{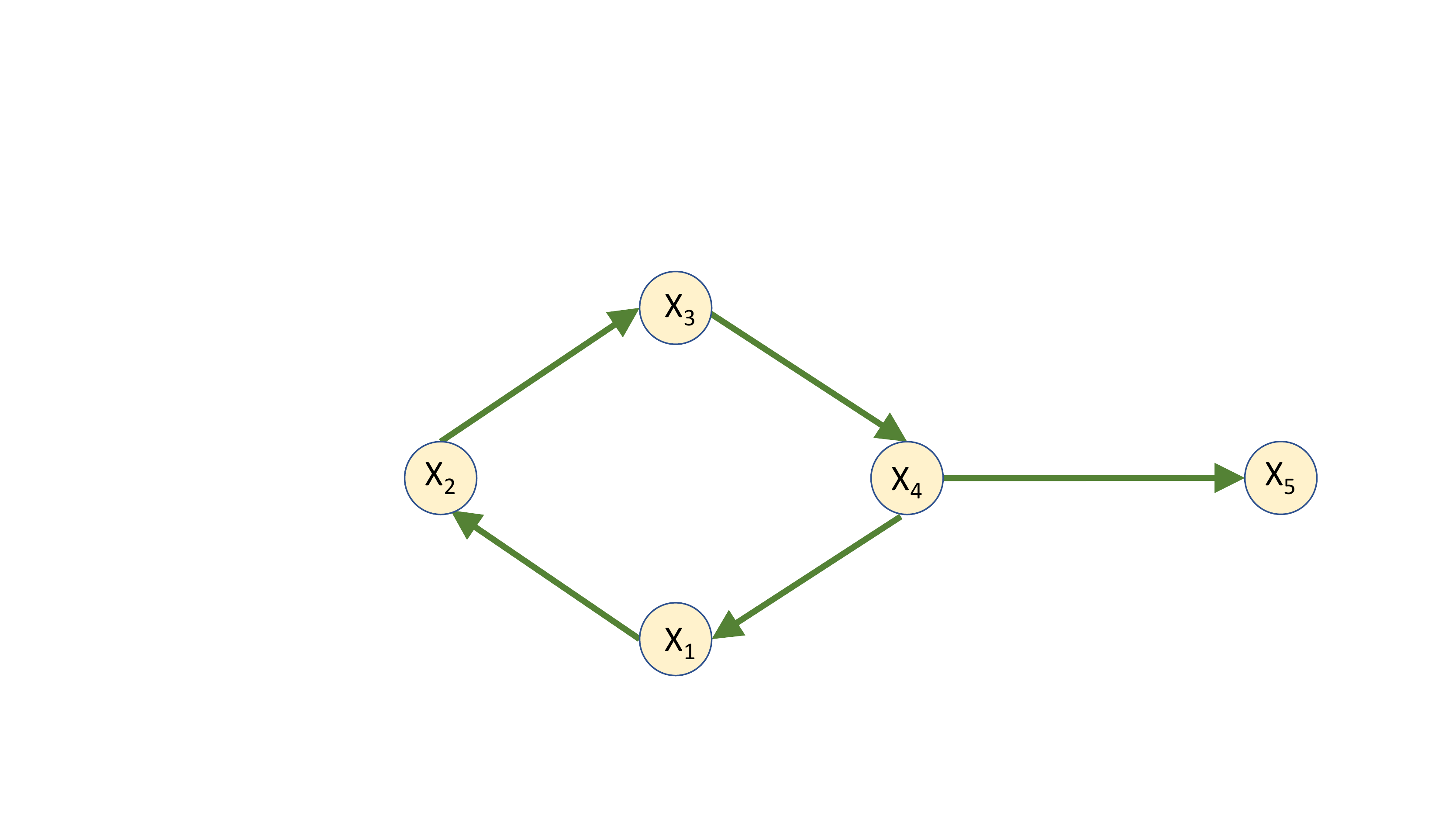}
\caption{Network corresponding to the dynamical system given in (\ref{feedback_eg}).}\label{feedback_network}
\end{figure}

{\bf Example 2}: In this example the objective is to identify network topology from time series data. Consider a network dynamical system described by following difference equation
{\small
\begin{eqnarray}\label{feedback_eg}
\begin{pmatrix}
z_{t+1}^1\\
z_{t+1}^2\\
z_{t+1}^3\\
z_{t+1}^4\\
z_{t+1}^5
\end{pmatrix} = 0.9\begin{pmatrix}
0 & 0 & 0 & 1 & 0\\
1 & 0 & 0 & 0 & 0\\
0 & 1 & 0 & 0 & 0\\
0 & 0 & 1 & 0 & 0\\
0 & 0 & 0 & 1 & 0
\end{pmatrix}\begin{pmatrix}
z_t^1\\
z_t^2\\
z_t^3\\
z_t^4\\
z_t^5
\end{pmatrix}+\sigma \xi_t
\end{eqnarray}} 
The network topology corresponding to above system is shown in Fig. \ref{feedback_network} 

Following Remark \ref{remark_topology}, we know that the information transfer from $T_{x_1\to y}$ is zero if and only if $A_{yx_1}=0$. Hence, information transfer can be used to inference presence or absence of a link i.e., network topology. In the above system $z^i$ dynamics affect the dynamics of $z^{i+1}$ for $i = 1,2,3,4$ and $z^1$ is affected by $z^3$. Hence, $T_{z^i\to z^{i+1}}$ should be non-zero for $i=1,2,3,4$ and $T_{z^3\to z^1}$ should also be non-zero. All the other transfers should be zero. The data was generated by choosing a single initial condition and propagating it for ten time steps with value of noise variance $\sigma=2.1$.

\begin{figure}[htp!]
\centering
\includegraphics[scale=.35]{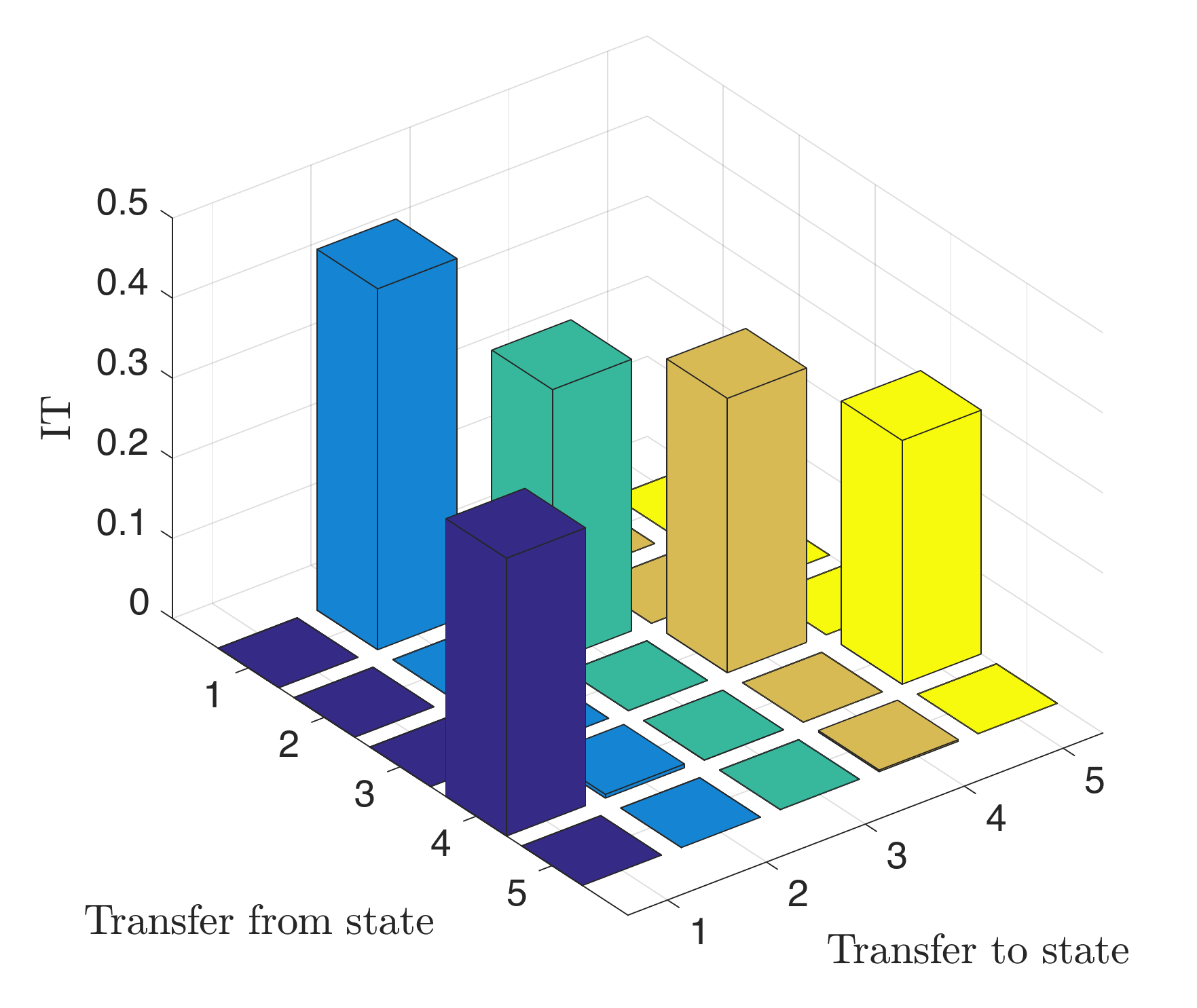}
\caption{Information transfer between the states.}\label{IT_fig_feedback}
\end{figure}

\begin{table}[htp!]
\centering
\caption{Information Transfer between the States}\label{IT_table}
\begin{tabular}{|c|c|c|c|}
\hline
I.T. & Value\\
\hline
$T_{z_1\to z_2}$ & 0.45\\
$T_{z_2\to z_3}$ & 0.34\\
$T_{z_3\to z_4}$ & 0.34\\
$T_{z_4\to z_5}$ & 0.31\\
$T_{z_4\to z_1}$ & 0.35\\
Other transfers & $\sim 10^{-4}$\\
\hline
\end{tabular}
\end{table}

The information transfer between the states is shown in table \ref{IT_table} and Fig. \ref{IT_fig_feedback} and we find that there is non-zero information transfer from $z^i\to z^{i+1}$ for $i=1,2,3,4$ and $T_{z^4\to z^1}$ is also non-zero. All the other transfers are very close to zero. For example, the transfer from $z^5$ to all the other states is of the order of $10^{-4}$ and hence we conclude that $z^5$ is not affecting any other state. Hence, we find that our information transfer measure recovers the correct causal structure or the network structure of the dynamical system from time series data. 

\begin{figure}
\centering
\includegraphics[scale=.3]{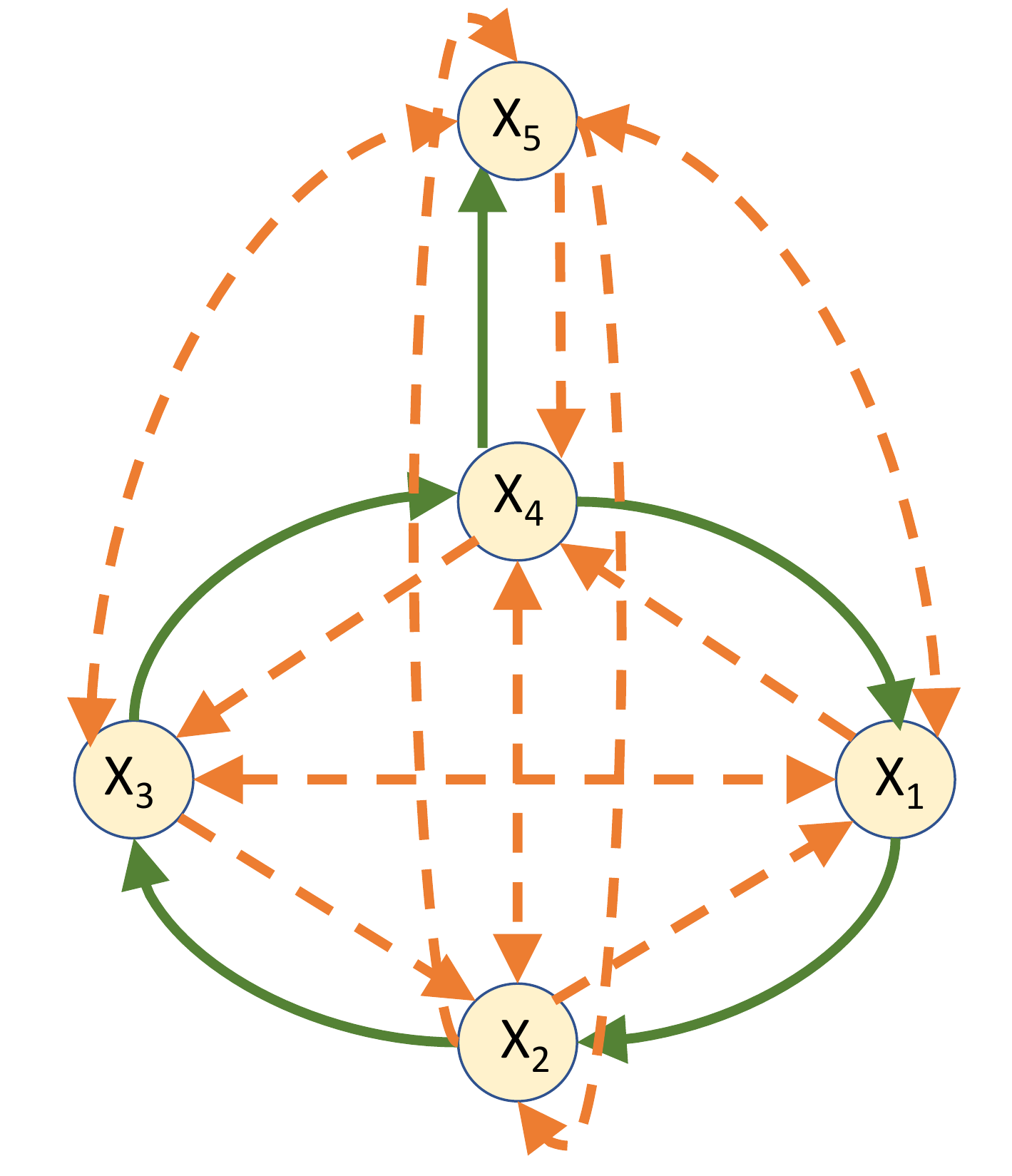}
\caption{Network identified by Granger causality}\label{granger_network}
\end{figure}

Further, to compare with an existing measure of causality, we used Granger causality test \cite{granger_causality},\cite{granger1969} on the same data set. Granger causality test is one of the most commonly used methods for causality detection and the intuition behind the definition is the following. A variable $X$ \emph{Granger causes} another variable $Y$ if the prediction of $Y$ based on its own past and past of $X$ is better than the prediction of $Y$ based on its own past alone. For self containment of the paper, we discuss the formulation of Granger causality briefly. Suppose $X_t$, $Y_t$ and $Z_t$ are three jointly distributed stationary multivariate stochastic processes. Consider the regression models 
\begin{eqnarray}\label{granger_regression}
X_t &=& \alpha_t + (X_{t-1}^{(p)}\oplus Z_{t-1}^{(r)})\cdot A + \epsilon_t
\end{eqnarray}
\begin{eqnarray}\label{granger_regression1}
X_t &=& \alpha_t^{\prime} + (X_{t-1}^{(p)}\oplus Y_{t-1}^{(q)}\oplus Z_{t-1}^{(r)})\cdot A^{\prime} + \epsilon_t^{\prime}
\end{eqnarray}
where $A$ and $A^{\prime}$ are the regression coefficients, $\alpha$ and $\alpha^{\prime}$ are constant terms, $\epsilon$ and $\epsilon^{\prime}$ are residuals, and the predictee variable $X$ is regressed first on the previous $p$ lags of itself plus $r$ lags of the conditioning variable $Z$ and second, in addition, on $q$ lags of the predictor variable $Y$. Granger causality of $Y$ to $X$, given $Z$, is a measure of of the extent to which the inclusion of $Y$ in the model (\ref{granger_regression1}) reduces the prediction error of the first model (\ref{granger_regression}) and is defined as
\begin{eqnarray}
G_{Y\to X|Z} = \ln\frac{\textnormal{var}(\epsilon_t)}{\textnormal{var}(\epsilon_t^{\prime})}
\end{eqnarray}
where $\textnormal{var}(\cdot)$ is the variance.

\begin{figure}[htp!]
\centering
\includegraphics[scale=.3]{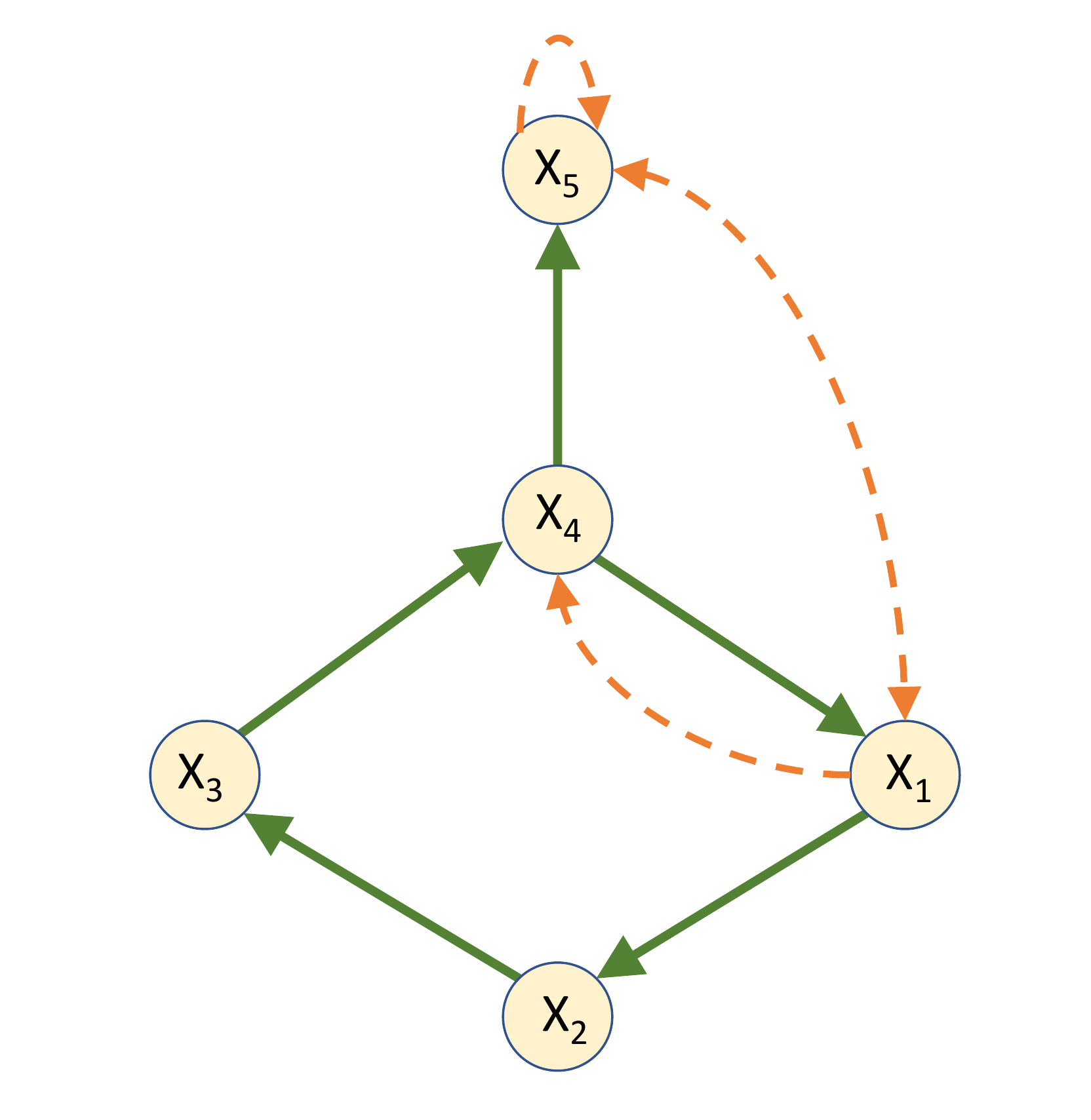}
\caption{Network obtained using Sparse DMD}\label{sparse_DMD}
\end{figure}

As can be seen from Fig. \ref{granger_network}, Granger causality identifies both direct and indirect causal influence and it fails to differentiate between direct and indirect influence. The direct links are shown in green and the indirect links, which are also identified by Granger causality are indicated by dotted orange lines. Hence, it is not possible to infer the correct causal structure from Granger causality test. However, our measure of information transfer can differentiate between direct and indirect influence \cite{sinha_IT_CDC2016} and information transfer computed in this paper gives only the direct influence. Hence the proposed measure captures the true causal structure.

At this stage one might wonder as to the need of computing information transfer for the purpose of identifying network topology. If all we are interested in determining the presence of absence of links between nodes, can one do that by simply estimating the system dynamics matrix $\bar A$ using optimization formulation in Eq. (\ref{regular2}) from time-series data. Is there a need to estimate $\bar A_{\not x}$ for the freeze system ?. In particular, if $\bar A_{ij}=0(\neq 0)$ then it implies absence (presence) of directed link from node $j$ to node $i$. To verify this claim in Fig. \ref{sparse_DMD}, we compare the results for the network topology identification obtained using our information transfer based method and one based on estimated system matrix $\bar A$. We find that the information transfer measure can regenerate the correct topology of the network, whereas, sparse DMD algorithm identifies some links which are not there in the original network. The spurious links are marked by dotted orange lines in Fig. \ref{sparse_DMD}.





{\bf Example 3}: Small world networks are ubiquitous in nature and in this example, we look at a small world network with 20 nodes and analyze how information transfer measure performs to recover the causal structure of the small world network. In the previous example, we constructed a network from a given dynamical system. In this example, we start with a network and construct a dynamical system from the network and identify the connections of the network using the information transfer measure. In particular, given a network of $n$ nodes, we construct a $n$-dimensional discrete time linear dynamical system such that, if there is a directed edge from node $i$ to node $j$, then $A_{ji}\neq 0$. Thus we obtain the system matrix of the dynamical system and we make the system stable by scaling the $A$ matrix with an appropriate constant.

\begin{figure}[htp!]
\centering
\subfigure[]{\includegraphics[scale=.17]{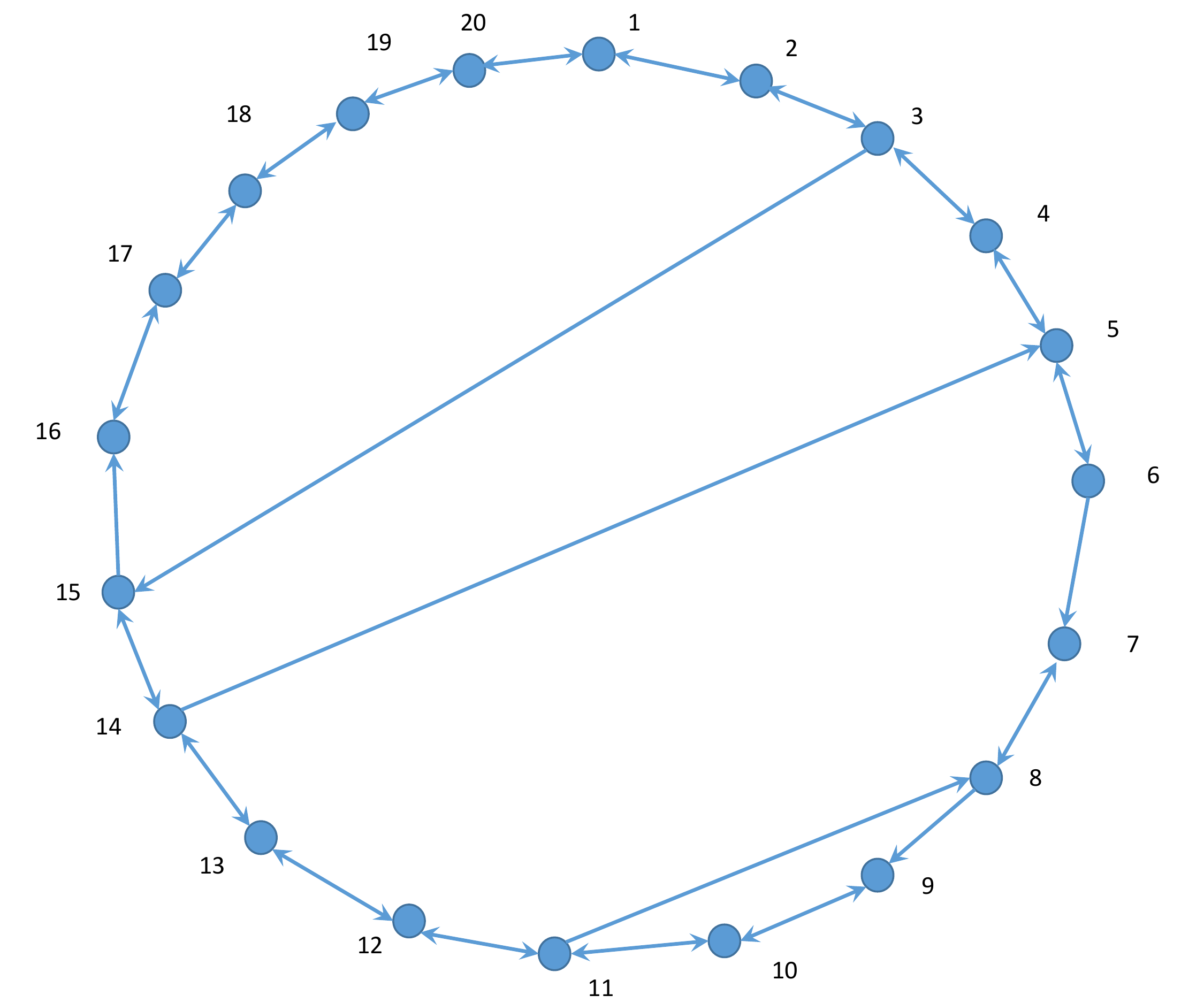}}
\subfigure[]{\includegraphics[scale=.17]{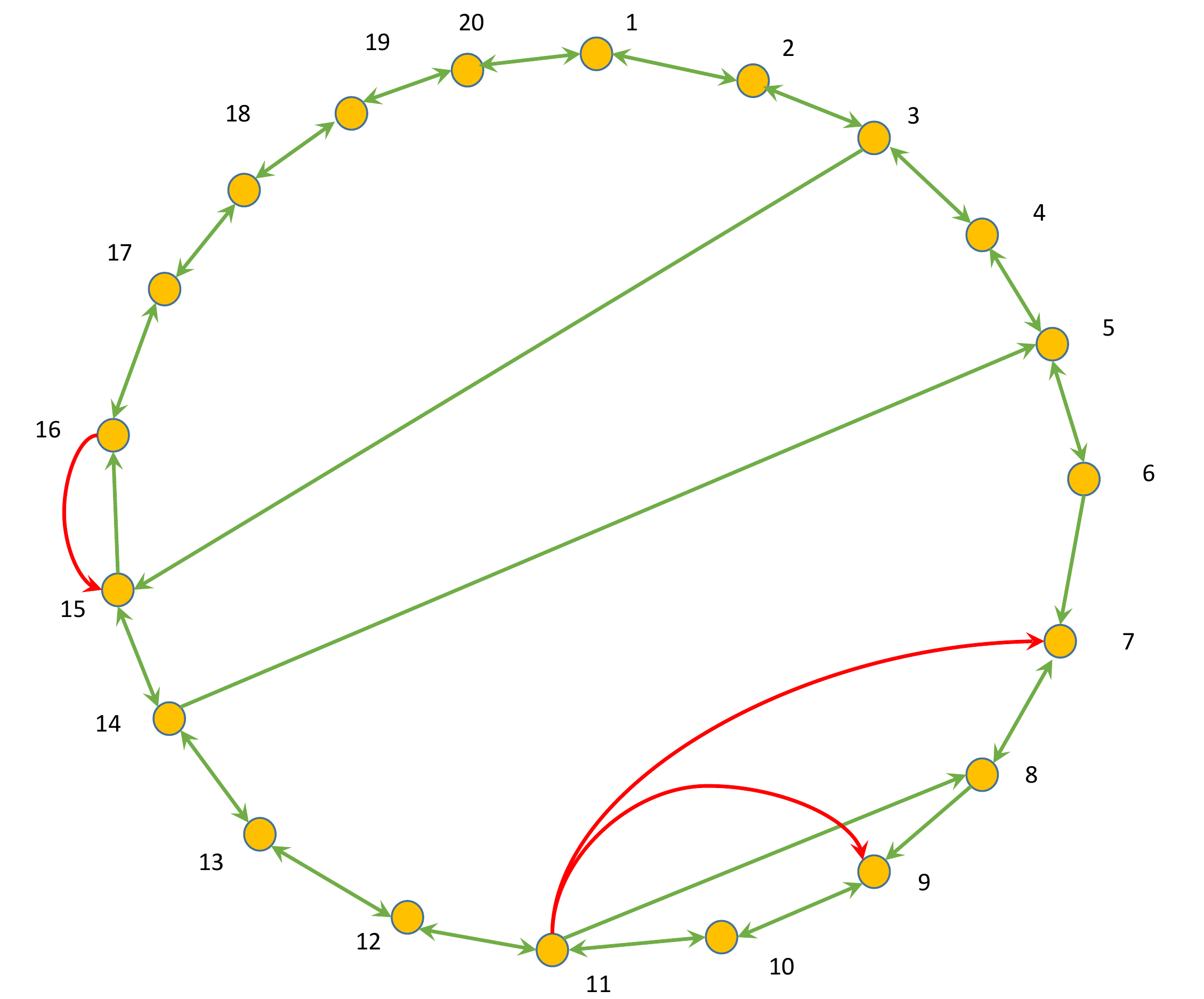}}
\caption{(a) Small world netowrk of 20 nodes. (b) Reconstructed network.}\label{20_node_small_world}
\end{figure}

In Fig. \ref{20_node_small_world}(a) we show the small world network of 20 nodes which is used to generate the data. In this case, we corrupt the data with i.i.d. Gaussian noise of variance $0.1$. Hence, $\lambda = 0.3$ in the optimization problem (\ref{regular}). Fig. \ref{20_node_small_world}(b) shows the reconstructed network using information transfer method. We find that in this case information transfer identifies three extra links, which are not there in the original system. These links are marked in red in Fig. \ref{20_node_small_world}(b). In this case, information transfer does not identify the exact causal structure and there is some error, but it performs pretty well.

\begin{figure}[htp!]
\centering
\includegraphics[scale=.25]{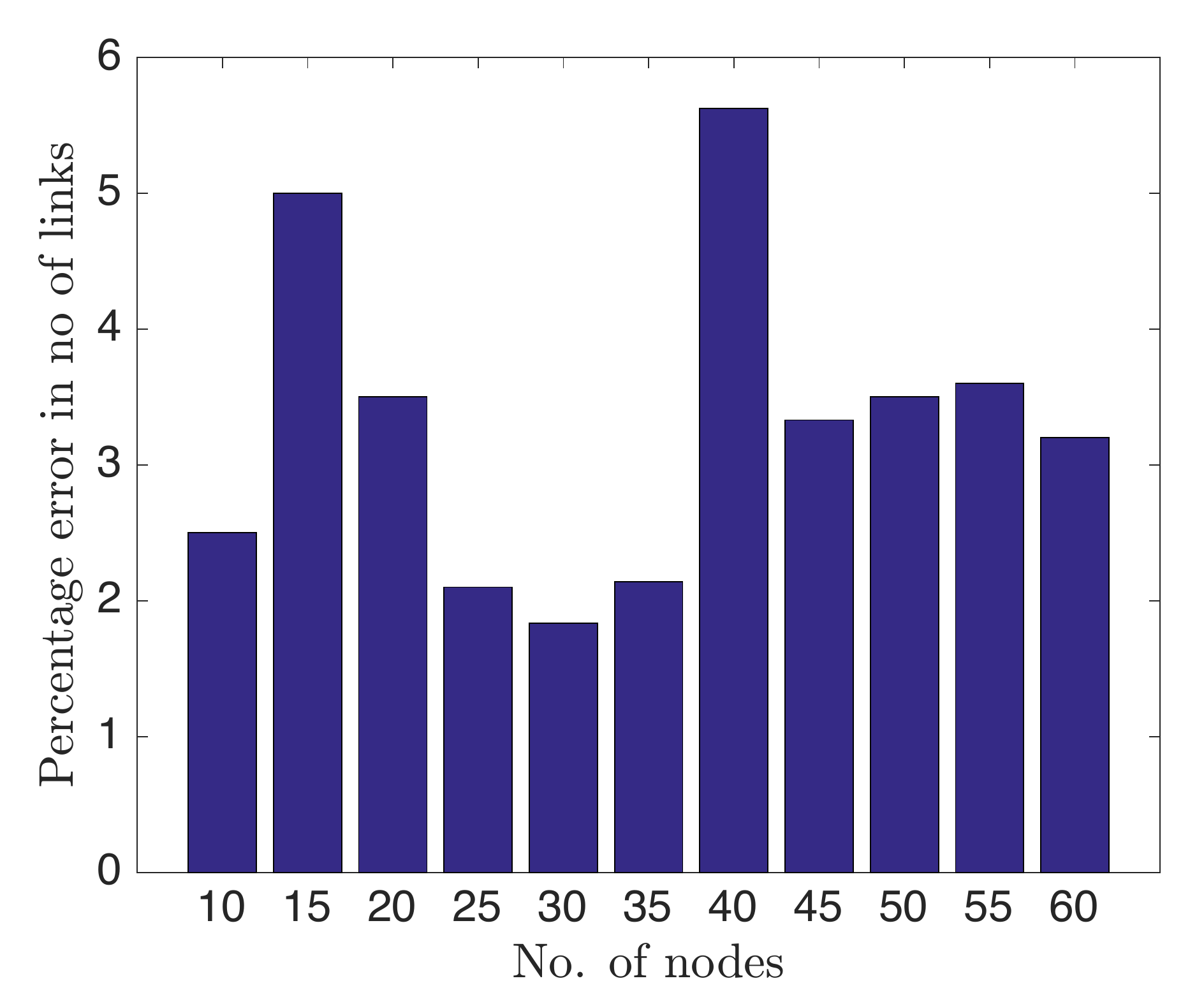}
\caption{Percentage error in number of links v/s number of nodes.}\label{perc_err_varying_nodes}
\end{figure}

In Fig. \ref{perc_err_varying_nodes}, we show the percentage error as the number of nodes in the system is increased. In all these cases, we did not corrupt the data with noise. In these examples, we find that even as the number of nodes is increased, the percentage error remains relatively small, and hence we conclude that the information transfer measure does a pretty good job in identifying the causal structure in a dynamical system. In these simulations, we considered small world networks of different number of nodes and the number of links is nearly four times the number of nodes. 
\section{Causal inference in nonlinear dynamical systems}\label{section_general}

One of the main challenge in computing the information transfer from time series data for nonlinear system is to propagation of probability density function under nonlinear flow field and the computation of conditional entropy term. Unlike linear system where linearity and Gaussian property of the probability density function was exploited for the data-driven computation of information transfer,  the same does not applies for nonlinear system. 
We make use of Naturally Structured Dynamic Mode Decomposition algorithm (NSDMD) for the approximation of transfer P-F operator and for propagation of probability density function. The positivity and the Markov property of the NSDMD algorithm is exploited for the propagation of probability density function and for computation of conditional entropy term.

Again we will outline the computation procedure for the two subspace case the general case (more than two subspace case) will follows from this procedure. The information transfer formula for the two subspace case can be simplified as follows. We rewrite Eq. (\ref{IT}) as
\begin{small}
\begin{eqnarray*}
&&T_{x\to y}=H(\rho(y(t+1)|y(t)))-H(\rho_{\not x}(y(t+1)|y(t)))\\
&& = H(\rho(y(t+1),y(t)))-H(\rho_{\not x}(y(t+1),y(t)))
\end{eqnarray*}
\end{small}
In writing the above equality we have used the fact that $H(X|Y) = H(X,Y)-H(Y)$. Furthermore, since $x$ is held frozen from time $t$ to time $t+1$, we have \[H(\rho(y(t)))=H(\rho_{\not{x}}(y(t)).\]

We next outline the procedure for computing the joint entropy term $H(\rho(y(t+1),y(t)))$ using the finite-dimensional approximation of the P-F matrix $\bf P$ obtained using NSDMD optimization formulation outlined in (\ref{cost}-\ref{positive3}). Note that in the construction of the P-F matrix $\bf P$ we use the original data set (\ref{data_original}). Once we outline the procedure for computing the joint entropy term for the non freeze case, $H(\rho(y(t+1),y(t)))$, the computation for the entropy term for the freeze case, $H(\rho_{\not x}(y(t+1),y(t)))$, will follow along similar lines. The only difference being the the P-F matrix for the freeze case, denoted by ${\bf P}_{\not x}$, will be computed using the modified data set obtained to replicate the freeze case i.e., data-set (\ref{mod_data}).

For computing $H(\rho(y(t+1),y(t)))$, we first consider finite approximation of $\rho(z(t+1),z(t))$ as discrete probability measure. Towards this goal we consider finite partition of the state space $Z$ as 
\[{\cal Z}=\{D_1,\ldots,D_K\}, \;\;\;{\cal D}=\cup_{k=1}^K D_k.\]
such that $D_i\cap D_j=\emptyset$.
Similarly, let
\[{\cal Z}^x=\{D_1^x,\ldots,D_K^x\},\;\;{\cal Z}^y=\{D_1^y,\ldots,D_K^y\}\]
where ${\cal Z}^x$ and ${\cal Z}^y$ are the projection of the partition, ${\cal Z}$, along the $x$ and $y$ coordinates respectively. Let ${\cal D}^x=\cup_{k=1}^K D_k^x$ and ${\cal D}^y=\cup_{k=1}^K D_k^y$. 
We have $\rho(z(t+1),z(t))=\rho(z(t))\rho(z(t+1)|z(t))$, Let 
\[[p_z]^t_{ij}:={\rm Prob}(z_{t+1}\in D_j|z_t\in D_i)\]
\[[p_z]_i^t:={\rm Prob}(z_t\in D_i)\]
\[{\rm Hence},\;\;{\rm Prob}(z_{t+1}\in D_j,z_t\in D_i)=[p_z]_i^t[p_z]^t_{ij}\]
Similarly, we can define $[p_x]_{i}^t$, $[p_y]_i^t, [p_x]_{ij}^t$, and $[p_y]_{ij}^t$ as follows:
\[[p_x]_i^t:={\rm Prob}(x_t\in D_i^x), \;\;[p_y]_i^t:={\rm Prob}(y_t\in D_i^y)
\]
To compute the above defined quantities, we make use of finite dimensional approximation of P-F operator, $\bf P$. Note that in the finite dimensional approximation of the P-F matrix using NSDMD algorithm we have assumed that the dictionary functions are positive i.e., $\psi_i(z)\geq 0$ for $i=1,\ldots K$. Furthermore, we also assume that the dictionary functions are density functions i.e.,
\[\int_Z \psi_i(z)dz=1,\;\;i=1,\ldots,K.\]
Let $w_t=(w_t^1,\ldots,w_t^K)\in \mathbb{R}^K$ be a probability row vector. Density function, $\rho(z(t))$ can be constructed using this probability vector and the dictionary functions ${\bf \Psi}(z)$ as $\rho(z(t))=w_t{\bf \Psi}^\top(z)$. This density function can be propagated using finite dimensional $\bf P$ as follows.
\[\rho(z(t+1))=w_t{\bf P}^\top {\bf \Psi}^\top(x)=w_{t+1}{\bf \Psi}^\top(x).\]

Hence we have 
\[[p_z]_i^t=\int_{D_i}w_t{\bf \Psi}^\top(z)dz=w_t \int_{D_i}{\bf \Psi}^\top(z)dz
= w_t \Theta_i^\top\]
where $\Theta_i = \int_{D_i}{\bf \Psi}(z)dz$.

Let $\lambda_i=\{i_1,i_2,\cdots , i_L\}\subseteq\{1,2,\cdots , K\}$ such that $D_{i_k}^y\cap D_i^y\neq \phi$, where $\phi$ is the empty set. Then
\begin{small}
\begin{eqnarray*}
[p_y]_i^t=\int_{D_{\lambda_i}}w_t{\bf \Psi}^\top (z)dz = w_t\int_{D_{\lambda_i}}{\bf \Psi}^\top (z)dz.
\end{eqnarray*}
\end{small}
Now let $\bar w$ be the probability vector such that 
\[\bar w\int_{D_i} {\bf \Psi}^\top(z)dz=1.\]
i.e., density function of the form $\bar w {\bf \Psi}^\top (z)$ correspond to the case where the entire distribution is concentrated on set $D_i$. With the above definition of $\bar w$ we have 
\[[p_z]_{ij}^t=\int_{D_j}\bar w{\bf P}^\top{\bf \Psi}dz=\bar w{\bf P}^\top\int_{D_j}{\bf \Psi}(z)dz = \bar{w}{\bf P}^\top\Theta_j^\top\]
 
Hence, we have 
\begin{eqnarray}\label{joint_prob}\nonumber
&&{\rm Prob}(z_{t+1}\in D_j,z_t\in D_i)]\\ \nonumber
&=&\left[w_t \int_{D_i}{\bf \Psi}^\top(z)dz\right] \left[\bar w{\bf P}^\top\int_{D_j}{\bf \Psi}(z)dz\right]\\
&=& w_t \Theta_i^\top \bar{w}{\bf P}^\top \Theta_j^\top  = \Gamma_{ij}
\end{eqnarray}

As defined earlier, the set $\mathcal{D}=\{\psi_1,\cdots , \psi_K\}$ are the dictionary functions for observables on the space $Z$. Note that $\rho(z(t+1),z(t))$ is defined on the product space $Z\times Z$ and hence we consider the set ${\Phi} = \mathcal{D}\times\mathcal{D}=\{\varphi_{11},\varphi_{12},\cdots ,\varphi_{KK}\}$ as the set of dictionary functions on the product space, where $\varphi_{ij} = \psi_i\psi_j$ and let $D_i\times D_j:=D_{i,j}$. Hence, we have 
\begin{eqnarray}
\rho(z(t+1),z(t))=\sum_{i,j=1}^K\Gamma_{ij}\varphi_{ij}(z,w).
\end{eqnarray}
 Let 
\begin{small}
\begin{eqnarray*}
\lambda_{i,j}=\{(p,k)|D_{i,j}^y\cap D_{p,k}^y\neq \phi ; p,k=1,2,\cdots , K\}
\end{eqnarray*}
\end{small}
Hence, by similar arguments for finding the marginal probability,
\begin{small}
\begin{eqnarray}\nonumber\label{conditional_y}
&&\textnormal{Prob}(y(t+1)\in D_j,y(t)\in D_i) \nonumber\\
&=& \int_{D_{\lambda_{i,j}}}\sum_{i,j = 1}^K\Gamma_{ij}\varphi_{ij}(z,w)dzdw
\end{eqnarray}
\end{small}
 Hence, (\ref{conditional_y}) gives the probability distribution of $(y(t+1),y(t))$ and using the entropy formula for a discrete probability distribution $(H(P)=-\sum_iP_i\log P_i)$, we get the entropy of $(y(t+1),y(t))$. Similarly, using the same above procedure for the modified data set (\ref{mod_data}), we can compute the entropy of $(y(t+1),y(t))$, when $x$ is held frozen and computing the difference $H(y(t+1),y(t))-H_{\not{x}}(y(t+1),y(t))$, we get the information transfer from $x$ to $y$.  

\begin{small}
\begin{algorithm}\label{algo2}
\caption{Algorithm for finding the information transfer from time series data}\label{algorithm}
\begin{enumerate}
\item{Compute the Koopman operator from the time series data using the method of Naturally Structured Dynamic Mode Decomposition.}
\item{Compute the joint probability of $(z(t+1),z(t))$ using (\ref{joint_prob}).}
\item{Compute the marginal probability of $(y(t+1),y(t))$ from the joint density from (\ref{conditional_y}).}
\item{Compute the entropy $H(y(t+1),y(t))$.}
\item{Form the modified data set from the given time series data from (\ref{mod_data}).}
\item{Repeat steps (1)-(5) for the modified data set to get the entropy $H_{\not{x}}(y(t+1),y(t))$.}
\item{Compute $H(y(t+1),y(t))-H_{\not{x}}(y(t+1),y(t))$ to get the information transfer from $x$ to $y$.}
\end{enumerate}

\end{algorithm}
\end{small}


\subsection{Examples and Simulations}\label{section_examples}

\begin{example} \emph{Two State Non-linear System.}
\newline
Next we consider a non-linear example. Consider the system 
\begin{eqnarray*}
x_{t+1} = 2x_t(1-x_t)+2y_t; \quad y_{t+1} = .8y_t
\end{eqnarray*}

The system was evolved for 300 time steps, starting from $\begin{pmatrix}
0.9 & 0.9
\end{pmatrix}^\top$. 

 From the system equations, we see that the $y$ dynamics is not affected by $x$, whereas, $x$ dynamics is affected by $y$. So there should be non-zero flow of information from $y$ to $x$ and there should be zero information flow from $x$ to $y$. 

For this example too, we considered Gaussian radial basis functions, with $\sigma=0.01$, for computation of the Koopman operators and with our algorithm we found $T_{x\to y}= -0.03$ and $T_{y\to x}=1.63$. So, in this case also our information transfer measure identifies that $y$ affects $x$ dynamics and $x$ does not influence $y$. Hence we have identified the correct causal structure. 

Granger causality for this example identified a statistical dependence of $x$ and $y$ and inferred that both $x$ and $y$ cause each other, whereas, in reality the $y$ dynamics is never affected by $x$ dynamics and hence the influence of $x$ on $y$ should be zero. However, Granger causality fails to identify this, while our information transfer measure does capture the zero influence. 

\end{example}

\begin{example} \emph{Henon Map.} In the second example, we consider the Henon map. It is one of the most studied dynamical systems which exhibit chaotic behaviour. The dynamical equations of the Henon map are 
\begin{eqnarray*}
x_{n+1} &=& 1 - ax_n^2 + y_n + \gamma\xi_x\\
y_{n+1} &=& bx_n + \gamma\xi_y
\end{eqnarray*}

The classical values of the parameters are $a = 1.4$ and $b = 0.3$. We add small process noise in the system with $\gamma = 0.01$. The support of the the first two Koopman eigenfunctions, computed using the NSDMD algorithm, is shown in Fig. \ref{henon_fig}.

\begin{figure}[htp!]
\centering
\subfigure[]{\includegraphics[scale = .19]{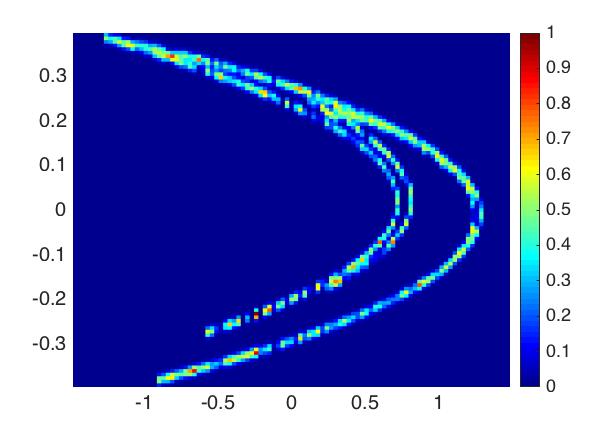}}
\subfigure[]{\includegraphics[scale = .19]{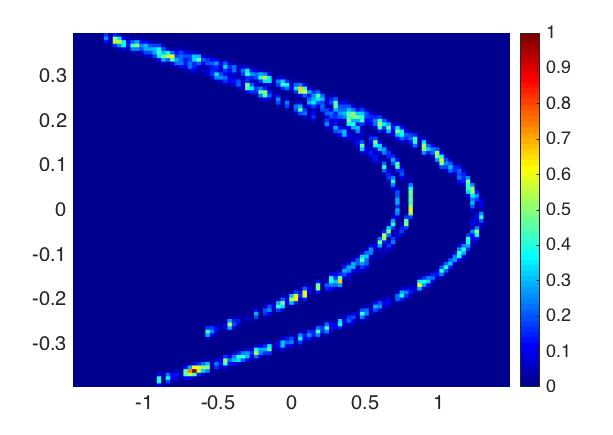}}
\caption{(a) Support of the eigenfunction corresponding to the largest eigenvalue of the Koopman operator. (b) Support of the eigenfunction corresponding to the second largest eigenvalue of the Koopman operator.}\label{henon_fig}
\end{figure}

For computing the information transfers, data was collected for 1000 time steps and information transfer between the states was calculated using the algorithm described in algorithm 2. We chose 200 Gaussian radial basis functions, with $\sigma = 0.01$ as the dictionary functions and this gave $T_{x\to y} = 0.0982$ and $T_{y\to x} = - 0.7246$. It is to be noted that one of the transfers is negative. We suspect that the negative value of the transfer is due to the fact that one Lyapunov exponent of the Henon map is negative, while the other is positive. Another interesting observation is the fact that the information transfer from $x$ to $y$ is nearly $b^2$. This is observed to be true in linear systems which have the following structure
\begin{eqnarray*}
x_{n+1} &=& a_x x_n + a_{xy} y_n + \sigma \xi_x\\
y_{n+1} &=& a_{yx}x_n + \sigma \xi_y
\end{eqnarray*}
In particular, it was observed that with small $\sigma \approx 0.01$, the information transfer from $x$ to $y$ is nearly equal to $a_{yx}^2$. 
Since the $x$ dynamics affect the $y$ dynamics linearly in the case of a Henon map, the same result hold true here as well. 


\end{example}

\section{Conclusion}\label{section_conclusion}
In this paper we address the problem of causal inference from time series data. Causality and influence characterization is an important problem and the existing measures of causality often fail to capture the true causal structure of a dynamical system. Based on a new definition of causality, which has been shown to capture the true causal structure, in this paper, we have provided a novel approach to identify the causal structure in a dynamical system. The general method is based on operator theoretic techniques for data analysis and requires the computation of Koopman operator via a newly developed scheme called Naturally Structured Dynamic Mode Decomposition. We have provided a complete algorithm to calculate information transfer in a dynamical system and use it to infer causality. We demonstrate our method on two different examples and show that this method does capture the true causal structure. Further, we provide a separate method to infer causality in linear systems and have shown how our measure of information transfer recovers the correct causal structure in a dynamical system.

\bibliographystyle{IEEEtran}
\bibliography{subhrajit_ref,ref,ref1}

\begin{thebibliography}{10}
\providecommand{\url}[1]{#1}
\csname url@samestyle\endcsname
\providecommand{\newblock}{\relax}
\providecommand{\bibinfo}[2]{#2}
\providecommand{\BIBentrySTDinterwordspacing}{\spaceskip=0pt\relax}
\providecommand{\BIBentryALTinterwordstretchfactor}{4}
\providecommand{\BIBentryALTinterwordspacing}{\spaceskip=\fontdimen2\font plus
\BIBentryALTinterwordstretchfactor\fontdimen3\font minus
  \fontdimen4\font\relax}
\providecommand{\BIBforeignlanguage}[2]{{%
\expandafter\ifx\csname l@#1\endcsname\relax
\typeout{** WARNING: IEEEtran.bst: No hyphenation pattern has been}%
\typeout{** loaded for the language `#1'. Using the pattern for}%
\typeout{** the default language instead.}%
\else
\language=\csname l@#1\endcsname
\fi
#2}}
\providecommand{\BIBdecl}{\relax}
\BIBdecl

\bibitem{sinha_IT_CDC2016}
S.~Sinha and U.~Vaidya, ``Causality preserving information transfer measure for
  control dynamical system,'' \emph{IEEE COnference on Decision and Control},
  pp. 7329--7334, 2016.

\bibitem{sinha_IT_ICC}
------, ``On information transfer in discrete dynamical systems,'' \emph{Indian
  Control Conference}, pp. 303--308, 2017.

\bibitem{mooij2016distinguishing}
J.~M. Mooij, J.~Peters, D.~Janzing, J.~Zscheischler, and B.~Sch{\"o}lkopf,
  ``Distinguishing cause from effect using observational data: methods and
  benchmarks,'' \emph{Journal of Machine Learning Research}, vol.~17, no.~32,
  pp. 1--102, 2016.

\bibitem{pearl_book}
J.~Pearl, \emph{Probabilistic reasoning in intelligent systems: networks of
  plausible inference}.\hskip 1em plus 0.5em minus 0.4em\relax Morgan Kaufmann,
  2014.

\bibitem{sprites_book_causation}
P.~Spirtes, C.~N. Glymour, and R.~Scheines, \emph{Causation, prediction, and
  search}.\hskip 1em plus 0.5em minus 0.4em\relax MIT press, 2000.

\bibitem{salapaka_reconstruction}
D.~Materassi and M.~V. Salapaka, ``On the problem of reconstructing an unknown
  topology via locality properties of the wiener filter,'' \emph{IEEE
  transactions on automatic control}, vol.~57, no.~7, pp. 1765--1777, 2012.

\bibitem{adebayo}
J.~Adebayo, T.~Southwick, V.~Chetty, E.~Yeung, Y.~Yuan, J.~Goncalves, J.~Grose,
  J.~Prince, G.-B. Stan, and S.~Warnick, ``Dynamical structure function
  identifiability conditions enabling signal structure reconstruction,'' in
  \emph{Decision and Control (CDC), 2012 IEEE 51st Annual Conference on}.\hskip
  1em plus 0.5em minus 0.4em\relax IEEE, 2012, pp. 4635--4641.

\bibitem{chow_sparse}
J.~Chow and P.~Kokotovic, ``Time scale modeling of sparse dynamic networks,''
  \emph{IEEE Transactions on Automatic Control}, vol.~30, no.~8, pp. 714--722,
  1985.

\bibitem{IT_socialmedia}
G.~V. Steeg and A.~Galstyan, ``Information transfer in social media,'' in
  \emph{{Proceedings of the 21st international conference on World Wide Web}},
  New York, NY, 2012, pp. 509--518.

\bibitem{IT_brain}
O.~Sporns, \emph{The networks of the brain}.\hskip 1em plus 0.5em minus
  0.4em\relax {MIT Press}, 2010.

\bibitem{IT_bionetwork1}
D.~J.~S. A.~Rao, A. O.~Hero and J.~D. Engel, ``Motif discovery in
  tissue-specific regulatory sequences using directed information.'' in
  \emph{{EURASIP J. on Bioinformatics and Systems Biology}}, 2007.

\bibitem{IT_bionetwork2}
------, ``Inference of biologically relevant gene influence networks using the
  directed information criterion.'' in \emph{{In proc. ICASSP}}, Toulouse,
  France, 2006.

\bibitem{granger_economics}
{C.W.J. Granger}, ``{Investigating causal relations by econometrics models and
  cross-spectral methods},'' \emph{{Econometrica}}, vol.~37, no.~3, pp.
  424--438, 1969.

\bibitem{IT_economics}
{C. A. Sims}, ``{Money, income and causality},'' \emph{{American Economic
  Review}}, vol.~62, pp. 540--552, 1972.

\bibitem{granger_causality}
C.~W.~J. Granger, ``Testing for causality,'' \emph{{Journal of Economic
  Dynamics and Control}}, vol.~2, pp. 329--352, 1980.

\bibitem{IT_massey_directed}
J.~L. Massey, ``Causality, feedback and directed information.'' in \emph{{Proc.
  Intl. Symp. on Info. th. and its Applications}}, Waikiki, Hawai, USA, 1990.

\bibitem{IT_kramer_directedit}
G.~Kramer, ``Directed information for channels with feedback,'' in \emph{{PhD
  Thesis}}, Swiss Federal Institute of Technology Zurich, 1998.

\bibitem{IT_schreiber}
T.~Schreiber, ``Measuring information transfer,'' \emph{{Physical Review
  Letters}}, vol. 85, no. 2, pp. 461--464, July, 2000.

\bibitem{tatonetti2012data}
N.~P. Tatonetti, P.~Y. Patrick, R.~Daneshjou, and R.~B. Altman, ``Data-driven
  prediction of drug effects and interactions,'' \emph{Science translational
  medicine}, vol.~4, no. 125, pp. 125ra31--125ra31, 2012.

\bibitem{kuhnert2014data}
C.~K{\"u}hnert and J.~Beyerer, ``Data-driven methods for the detection of
  causal structures in process technology,'' \emph{Machines}, vol.~2, no.~4,
  pp. 255--274, 2014.

\bibitem{shu2013data}
Y.~Shu and J.~Zhao, ``Data-driven causal inference based on a modified transfer
  entropy,'' \emph{Computers \& Chemical Engineering}, vol.~57, pp. 173--180,
  2013.

\bibitem{zou2009granger}
C.~Zou and J.~Feng, ``Granger causality vs. dynamic bayesian network inference:
  a comparative study,'' \emph{BMC bioinformatics}, vol.~10, no.~1, p. 122,
  2009.

\bibitem{liang_kleeman_prl}
{X. S. Liang and R. Kleeman}, ``{Information transfer between dynamical system
  components},'' \emph{{Physical Review Letters}}, vol.~95, p. 244101, 2005.

\bibitem{robust_dmd_acc}
S.~Sinha, B.~Huang, and U.~Vaidya, ``Robust approximation of koopman operator
  and prediction in ransdom dynamical systems,'' \emph{submitted to American
  Control Conference}, 2018.

\bibitem{umesh_website}
\BIBentryALTinterwordspacing
------, ``Robust approximation of koopman operator and prediction in ransdom
  dynamical systems.'' [Online]. Available:
  \url{http://home.eng.iastate.edu/~ugvaidya/publications.html}
\BIBentrySTDinterwordspacing

\bibitem{Lasota}
A.~Lasota and M.~C. Mackey, \emph{Chaos, Fractals, and Noise: Stochastic
  Aspects of Dynamics}.\hskip 1em plus 0.5em minus 0.4em\relax New York:
  Springer-Verlag, 1994.

\bibitem{williams2015data}
M.~O. Williams, I.~G. Kevrekidis, and C.~W. Rowley, ``A data--driven
  approximation of the koopman operator: Extending dynamic mode
  decomposition,'' \emph{Journal of Nonlinear Science}, vol.~25, no.~6, pp.
  1307--1346, 2015.

\bibitem{schmid2010dynamic}
P.~J. Schmid, ``Dynamic mode decomposition of numerical and experimental
  data,'' \emph{Journal of fluid mechanics}, vol. 656, pp. 5--28, 2010.

\bibitem{Umesh_NSDMD}
B.~Huang and U.~Vaidya, ``Data-driven approximation of transfer operators:
  Naturally structured dynamic mode decomposition,'' in
  \emph{{https://arxiv.org/abs/1709.06203}}, 2016.

\bibitem{caramanis201214}
C.~Caramanis, S.~Mannor, and H.~Xu, ``Robust optimization in machine
  learning,'' \emph{Optimization for machine learning}, p. 369, 2012.

\bibitem{granger1969}
C.~W. Granger, ``Investigating causal relations by econometric models and
  cross-spectral methods,'' \emph{Econometrica: Journal of the Econometric
  Society}, pp. 424--438, 1969.

\end{thebibliography}

\end{document}